\documentclass[11pt,letterpaper]{amsart}

\usepackage{amsmath,amssymb,amsthm,amsxtra,mathrsfs}
\usepackage{yhmath}
\usepackage{graphicx}
\usepackage{cancel,ulem}
\usepackage{color}
\usepackage[all,cmtip]{xy}
\usepackage{hyperref}
\usepackage{enumitem}
\usepackage{color}
\usepackage{bm}
  \usepackage{subcaption}
\usepackage{tikz}
\usetikzlibrary{arrows}
\usepackage{accents}
\usepackage{bbm}
\usetikzlibrary{arrows.meta, decorations.pathmorphing}
\usepackage{dsfont}

\newcommand{\la}{\langle}
\newcommand{\ra}{\rangle}

\newcommand{\jac}{\text{Jac}}
\newcommand{\vertiii}[1]{{\left\vert\kern-0.25ex\left\vert\kern-0.25ex\left\vert #1 
    \right\vert\kern-0.25ex\right\vert\kern-0.25ex\right\vert}}

\def\r{{\rangle}}
\def\l{{\langle}}
\def\S{{\mathbb S}}

\def\R{{\mathbb R}}

\def\K{\mathcal{K}}
\def\hat{\widehat}
\newtheorem{theorem}{Theorem}
\newtheorem{lemma}[theorem]{Lemma}

\newtheorem{proposition}[theorem]{Proposition}

\theoremstyle{definition}

\theoremstyle{remark}
\newtheorem{remark}[theorem]{Remark}

\numberwithin{equation}{section}
\numberwithin{theorem}{section}
\numberwithin{problem}{section}

\begin{document}

\begin{abstract}
We prove a new moment-preserving Young's inequality for the gain operator of the Boltzmann equation with hard potentials, including the critical case of hard-spheres. Our approach relies crucially on a novel cancellation mechanism dealing with the pathological case of energy-absorbing collisions (that is, collisions that accumulate energy to only one of the outgoing particles). This difficulty is specific to hard potentials, and is not present for Maxwell molecules. Our method quantifies the heuristic that, while energy-absorbing collisions occur with non-trivial probability, they are statistically rare, and therefore do not affect the overall averaging behavior of the gain operator. At the technical level, our proof relies solely on tools from kinetic theory, such as geometric identities and angular averaging.

\end{abstract}

\title[Moment-preserving Young's inequality for hard potentials]{Moment-preserving Young's inequality for the hard potential Boltzmann gain operator}

\author[I. Ampatzoglou]{Ioakeim Ampatzoglou}
     \address{Corresponding author: Ioakeim Ampatzoglou, Department of Mathematics, Baruch College \& The Graduate Center,  City University of New York, Newman Vertical Campus, 55 Lexington Ave New York, NY, 10010, USA}
 \email{ioakeim.ampatzoglou@baruch.cuny.edu ; iampatzoglou@gc.cuny.edu}

 \author[T. L\'eger ]{Tristan L\'eger}

\maketitle
\tableofcontents
\section{Introduction}

\subsection{Problem setup and goal of this paper}

The Boltzmann equation is the central equation of collisional kinetic theory and one of the most celebrated equations of mathematical physics. 
It describes dilute gases of microscopic interacting particles at the mesoscopic/kinetic level; more precisely, it predicts the time evolution of the probability density $f(t,x,v)$ of a dilute gas in non-equilibrium. The Boltzmann equation reads
\begin{equation}\label{classical Boltzmann equation}
\partial_t f+v\cdot\nabla_x f=Q(f,f),
\end{equation}
where the advection operator $\partial_t+v\cdot\nabla_x$ accounts for the rectilinear motion of particles, while $Q(f,f)$ is the quadratic collisional operator encoding the binary interactions between gas particles. It is given by
\begin{equation}\label{classical Boltzmann operator}
Q(f,f)(t,x,v)=\int_{\R^3\times\S^2}B(u,\sigma)\Big(f(t,x,v^*)f(t,x,v_1^*)-f(t,x,v)f(t,x,v_1)\Big)\,d\sigma\,dv_1,
\end{equation}
where $v,v_1$ denote the pre-collisional velocities of the incoming particles, $v^*$, $v_1^*$ denote the post-collisional velocities, and $u:=v-v_1$ denotes the relative velocity of the particles before the collision. Simplifying the notation slightly, we drop the dependence of the arguments on $t,x.$

In the case of perfectly elastic collisions, and assuming without loss of generality that all particles are of equal mass $m=1$, the post-collisional velocities are given by:
\begin{equation}\label{post-collisional velocities}
\begin{split}
v^*&=V-\frac{|u|}{2}\sigma,\qquad v_1^*=V+\frac{|u|}{2}\sigma,
\end{split}
\end{equation}
where $V:=\frac{v+v_1}{2}$ is the center of mass of the velocities of the incoming velocities and $\sigma\in\S^2$ represents the scattering direction of the collision. It is then straightforward to see that the momentum and energy are conserved, as well as the relative velocity magnitude:
\begin{align}
v^*+v_1^*&=v+v_1,  \label{conservation of momentum}\\
|v^*|^2+|v_1^*|^2&=|v|^2+|v_1|^2,\label{conservation of energy}\\
|v^*-v_1^*|&=|v-v_1|.\label{conservation of relative velocity}
\end{align}

The kernel $B(u,\sigma)$ appearing in \eqref{classical Boltzmann operator} is the collisional cross-section, which is assumed to be of the form 
\begin{equation}\label{cross-section}
    B(u,\sigma)=|u|^{\gamma}\, b(\widehat{u}\cdot\sigma),\quad -3<\gamma<2,
\end{equation}
where the angular cross-section $b:[-1,1]\to\R$ satisfies:
\begin{align}
b(z)&\geq 0,\quad\forall\, z\in[-1,1]\tag{H.1},\label{assumption 1}\\
b(-z)&=b(z),\quad \forall\,z\in[-1,1],\tag{H.2}\label{assumption 2}\\
&b\in L^\infty([-1,1]).\tag{H.3}\label{assumption 3}
\end{align}

The parameter $\gamma$ in the radial part of the cross-section is of crucial importance: it determines the type of intramolecular potential considered.

\begin{itemize}
    \item The case $0<\gamma< 2$ corresponds to hard potentials, where the interaction is strong and of short range. The physically meaningful range is $0<\gamma\leq 1$, which includes the celebrated hard-sphere model ($\gamma=1,\,b=1$), where particles interact like impenetrable billiard balls. The range $1<\gamma<2$ corresponds to super-hard potentials.  Hard-spheres are of particular importance since it was the model originally introduced by Boltzmann \cite{boltzmann} and Maxwell \cite{ma867}. It essentially remains the only model for which the Boltzmann has been rigorously derived from deterministic many particle dynamics \cite{Lanford,GSRT,DHM}. 
    \item The case $\gamma=0$ corresponds to Maxwell molecules. 
    \item The case $-3<\gamma<0$ corresponds to soft potentials. This type of interaction is smoother and of longer range.
\end{itemize}

Both physically and mathematically, these three types of potentials lead to different kinds of difficulties that typically require separate treatment.
In this paper we focus on hard potentials, i.e. we assume $0<\gamma<2$. In the course of the proof, we also improve existing results on Maxwell molecules with bounded cross-sections.

Regarding the assumptions on the angular cross-section $b$, 
condition \eqref{assumption 2} corresponds to indistinguishably of the particles, while the boundedness condition \eqref{assumption 3} is a slightly stronger assumption than the usual Grad's cut-off, which is the mere integrability condition
 $
 \int_{\S^2}b(\widehat{u}\cdot\sigma)\,d\sigma<\infty.
 $
 Yet \eqref{assumption 3} covers most cases of physical interest for hard potentials, and in particular the classical hard-sphere model ($\gamma=1,\,b=1$).
Most importantly, since \eqref{assumption 3} is a special case of Grad's cut-off, it implies that the collisional operator can be written in gain-loss form
\begin{equation} \label{splitting of Q}
Q(f,f)=Q^+(f,f)-Q^-(f,f),
\end{equation}
where the gain and loss operators are respectively given by 
\begin{align} 
Q^+ (f,f)&=\int_{\R^3\times\S^2}|u|^\gamma b(\widehat{u}\cdot\sigma) f(v^*)f(v_1^*)\,d\sigma\,dv_1,\label{gain operator}\\
Q^- (f,f)&=\int_{\R^3\times\S^2}|u|^\gamma b(\widehat{u}\cdot\sigma) f(v)f(v_1)\,d\sigma\,dv_1.\label{loss operator}
\end{align}

\subsubsection*{Half-sphere restriction and bilinear operators}
Exploiting spherical symmetry, it is often convenient to restrict the collisional operators to a half-sphere. For hard potentials, the most natural restriction is to fix the outgoing particle with the largest velocity (for Maxwell molecules we will restrict to the half-sphere corresponding to $\widehat{u}\cdot\sigma>0$ instead, see Section \ref{sec:maxwell}). Namely, by conservation of energy we have $|v^*|^2+|v_1^*|^2=|v|^2+|v_1|^2:=E$, thus either $|v_1^*|^2\geq E/2$ or $|v^*|^2\geq E/2$. Since the set where $|v_1^*|^2=|v^*|^2=E/2$ is of measure zero, we obtain
\begin{align}
Q^+(f,f)&=\int_{\R^3\times\S^2}|u|^\gamma b(\widehat{u}\cdot\sigma)f(v^*)f(v_1^*)\mathds{1}_{|v_1^*|^2>E/2}\,d\sigma\,dv_1\notag\\
&\hspace{1cm}+\int_{\R^3\times\S^2}|u|^\gamma b(\widehat{u}\cdot\sigma)f(v^*)f(v_1^*)\mathds{1}_{|v^*|^2>E/2}\,d\sigma\,dv_1\notag\\
&=\int_{\R^3\times\S^2}|u|^\gamma b(\widehat{u}\cdot\sigma)f(v^*)f(v_1^*)\mathds{1}_{|v_1^*|^2>E/2}\,d\sigma\,dv_1\notag\\
&\hspace{1cm}+\int_{\R^3\times\S^2}|u|^\gamma b(-\widehat{u}\cdot\sigma)f(v_1^*)f(v^*)\mathds{1}_{|v_1^*|^2>E/2}\,d\sigma\,dv_1\notag\\
&=2\int_{\R^3\times\S^2}|u|^\gamma b(\widehat{u}\cdot\sigma)f(v^*)f(v_1^*)\mathds{1}_{|v_1^*|^2>E/2}\,d\sigma\,dv_1,\notag
\end{align}
where in the second part of the summation we used the substitution $\sigma\to -\sigma$ (which just interchanges $v^*$ and $v_1^*$ by \eqref{post-collisional velocities}) and for the last equality we used \eqref{assumption 2}. Similarly, we obtain
\begin{equation}\label{quadratic loss restricred}
Q^{-}(f,f)=2\int_{\R^3\times\S^2}|u|^\gamma b(\widehat{u}\cdot\sigma)f(v)g(v_1)\mathds{1}_{|v_1^*|^2>E/2}\,d\sigma\,dv_1 .
\end{equation}

Motivated by this computation, we introduce the bilinear forms:
\begin{align}
Q(f,g)&:=Q^+(f,g)-Q^-(f,g),\label{full operator bilinear}\\
Q^+ (f,g)&:=2\int_{\R^3\times\S^2}|u|^\gamma b(\widehat{u}\cdot\sigma) f(v^*)g(v_1^*)\mathds{1}_{|v_1^*|^2>E/2}\,d\sigma\,dv_1,\label{gain operator bilinear}\\
Q^- (f,g)&:=2\int_{\R^3\times\S^2}|u|^\gamma b(\widehat{u}\cdot\sigma) f(v)g(v_1)\mathds{1}_{|v_1^*|^2>E/2}\,d\sigma\,dv_1.\label{loss operator bilinear}
\end{align}

Clearly, the loss operator $Q^-(f,g)$ is much simpler to analyze, and essentially behaves like the product $f(|v|^\gamma\ast g).$ 

In contrast, the gain operator $Q^+(f,g)$ has a much richer structure due to the averaging over the angular variable and exhibits regularization effects \cite{Gus2,Lions I}. As a result, it is the main object to understand when developing the mathematical theory of the Boltzmann and related equations.

Quantifying these regularizing effects is subtle, though, because of the complex averaging structure of the operator. The first step is to choose a suitable topology to measure the gain of regularity. The most relevant functional analytic setting is  polynomially weighted Lebesgue spaces only in the velocity variable $\l v \r^{-k} L^p_v,$  for several reasons:
\begin{itemize}
\item Such bounds are space translation invariant, thus they can be applied whether the spatial domain is bounded or unbounded. Indeed, $\l v \r^{-k} L^q_x L^p_v$ estimates easily follow from H\"{o}lder's inequality. This feature is important to address questions related to the hydrodynamic limit of \eqref{classical Boltzmann equation}. In that field, most mathematical results are proved in spatially bounded domains such as the torus. Without being exhaustive, we refer to the recent work of Deng-Hani-Ma \cite{DHM2} for a complete derivation of fundamental fluid dynamics equations from many-body classical dynamics. This followed their groundbreaking work \cite{DHM} on the long-time derivation of the Boltzmann equation from  hard-sphere dynamics which extended Lanford's theorem \cite{Lanford,GSRT} globally in time. 
\item Polynomial weights are the most suitable for studying the long-time behavior of the equation for large data, which is a major open problem, as outlined in \cite{Villani-review}. Large data results for the space homogeneous (i.e. space independent) Boltzmann equation for hard potentials have been proved by establishing generation and propagation of moments \cite{Pov,Elm,de93,we94,we97,bo97,miwe99}, which correspond exactly to Lebesgue spaces with polynomial weights in $v$. Moreover, as noted by L. Desvillettes and C. Villani \cite{DV}, the evolution cannot be treated by linearization techniques which typically require exponential weights. 
\item Working with polynomial weights is also more natural in the closely related field of wave turbulence, where equilibria are typically inverse polynomials in $v$ (e.g. Rayleigh-Jeans distributions or the celebrated Kolmogorov-Zakharov spectra, see e.g. Nazarenko \cite{Nazarenko} for more details). Therefore,  we expect the ideas from the present paper to also be helpful in studying  wave kinetic equations. They were already implemented our papers \cite{AmLe24} and \cite{AmLeNote}. The method proved quite versatile, allowing us to prove a global regularity result for small data, as well as an almost sharp local well posedness result for the kinetic wave equation. Both improve on existing literature (e.g. \cite{Ampatzoglou} and Germain-Ionescu-Tran \cite{GeIoTr}) while using the unified approach used in the present article.
\color{black}
\item Finally, the conditions for the solution to be derivable from hard-sphere dynamics, and thus solve Hilbert's sixth problem require  weights only in $v$, see the statement of Theorem 1 in the seminal work of Deng-Hani-Ma \cite{DHM}.
\end{itemize}

In this paper, we use the averaging present in the gain operator to prove new moment-preserving Young's inequalities in such polynomially weighted spaces. We work in the case of hard potentials, and the results include hard-spheres. Our approach relies crucially on a novel cancellation mechanism to deal with pathological collisions that accumulate energy to only one of the outgoing particles. As we shall see, this type of collisions distinguishes hard potentials from Maxwell molecules. Hence, a thorough understanding of these interactions is crucial for the study of hard potentials.

\subsubsection*{Related literature}
After the foundational works of Carleman \cite{Carleman 1, Carleman 2} ($1<p<\infty$)  and Arkeryd \cite{Arkeryd} ($p=\infty$) on the propagation of $L^p$--norms  of solutions to the homogeneous Boltzmann equation, it was Gustafsson \cite{Gus1,Gus2} who first noticed the convolution structure of the gain collisional operator by proving a Young-type inequality for $Q^+(f,g)$, using an intricate nonlinear interpolation procedure. Later, in a pioneering work, Lions \cite{Lions I} used Fourier integral operators techniques to show that $Q^+(f,g)$ gives a derivative gain. Wennberg \cite{we94reg} simplified the proof of Lions and arrived at the same estimates by using the Carleman representation and classical Fourier analysis tools. See also the works \cite{bode98,Lu98} for another family of estimates on the quadratic operator $Q^+(f,f)$  that are obtained by much simpler means. Based on these results, Toscani-Villani \cite{tovi00} and Mouhot-Villani \cite{MoVi04} studied the propagation of $L^p$--norms and smoothness of solutions to the homogeneous
Boltzmann equation. However, all these early results were obtained assuming point-wise truncations in the angular variable. More precisely, in these works singular collisions such as grazing $(\widehat{u}\cdot\sigma=0)$ or frontal $(\widehat{u}\cdot\sigma=\pm 1)$ are discarded.



More recently, Alonso-Carneiro \cite{AlCa10} used Fourier analysis and radial symmetrization techniques to obtain a Young-type inequality for Maxwell molecules and hard potentials, under integrability assumptions on the cross-section $b$ that replaced point-wise truncations. The dependence on the cross-section is embedded in the constants of their estimates. Special cases of those estimates were  previously obtained in \cite{DuKiRj06} and \cite{GaPaVi04}.  Soon after that, Alonso-Carneiro-Gamba \cite{AlCaGa10} extended the results of \cite{AlCa10} to inelastic collisions, and also proved a Hardy-Littlewood-Sobolev type of inequality for soft potentials ($-3<\gamma<0$).

However, in the case of hard potentials, the estimates obtained in \cite{AlCa10,AlCaGa10} are moment-increasing, due to the point-wise linear growth of the cross-section. 
More precisely for $\frac{1}{p}+\frac{1}{q}=1+\frac{1}{r}$, the estimate provided in \cite{AlCa10, AlCaGa10} is:
\begin{align*}
    \Vert \l v \r^k \,  Q^{+}(f,g) \Vert_{L^r_v} \lesssim \Vert \l v \r^{k+\gamma} \,  f \Vert_{L^p_v} \Vert \l v \r^{k+\gamma} \, g \Vert_{L^q_v}, \quad \l v\r:=\sqrt{1+|v|^2},
\end{align*}
which is not sufficient to prove local/global well-posedness of the Boltzmann equation in spaces of the form $\l v \r^{-k} L^q_x L^p_v$, or propagation of moments. For this, one needs moment-preserving estimates, and the main difficulty in obtaining those is precisely this point-wise linear growth of the potential.


With the collisional averaging mechanism  we uncover in this paper, we are able to offset this linear growth  and produce a new wide class of polynomially weighted moment-preserving Young's inequalities.

\subsection{Heuristics}\label{subsec:heuristics}
For simplicity, in this subsection we focus on the case of hard-spheres i.e. $\gamma=1$ and $b=1$, which is the prototypical hard potential model. We write the weighted \eqref{gain operator bilinear} as 
\begin{align*}
 \l v\r^k Q^+(f,g)(v) = 2 \int_{\mathbb{R}^3 \times \S^2} \frac{\vert u \vert \l v \r^k}{\l v^* \r^k \l v_1^* \r^k} f_k(v^*) g_k(v_1^*) \mathds{1}_{|v_1^*|^2>E/2
 } \, d\sigma dv_1 .
\end{align*}
The central object in our analysis is the ratio 
\begin{equation}\label{ratio intro}
\mathcal{R}:=\frac{|u| \l v\r^k}{\l v^*\r^k \l v_1^*\r^k},\quad u=v-v_1,\quad (v,v_1,\sigma)\in \R^3\times\R^3\times\S^2.
\end{equation}
We distinguish two types of collisions: \textit{energy-distributing} and \textit{energy-absorbing}. Each of these requires a different argument.

\subsubsection*{Energy-distributing collisions} 
\begin{figure}[h!]
    \centering
    \includegraphics[width=0.3\textwidth]{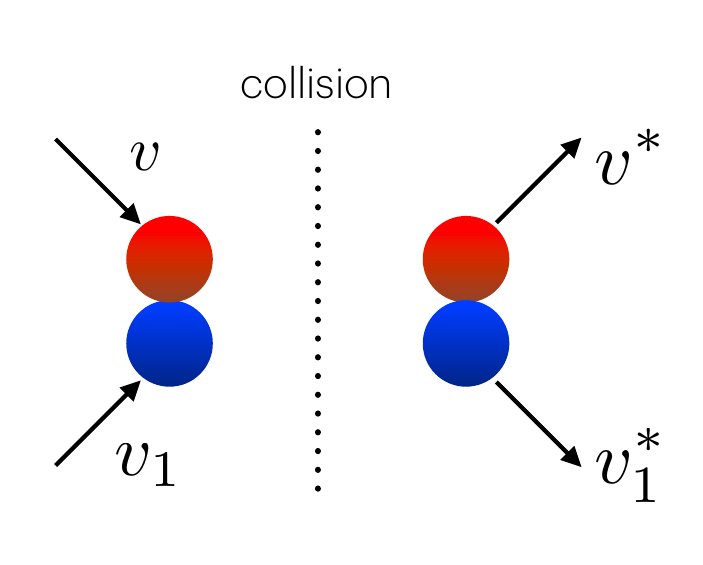} 
    \caption{Energy-distributing collisions i.e. $|v^*|\approx|v_1^*|\approx E^{1/2}$}
    \label{energy distributing figure}
\end{figure}

In most cases, energy is distributed among the outgoing particles, i.e.  $|v^*|^2$ is comparable to  the energy  as well. We call these collisions \textit{energy-distributing}. In Figure \ref{energy distributing figure}, we show an example of such a collision.

Now, since $|u|\leq |v|+|v_1|\lesssim E^{1/2}$, we can easily estimate
$$\mathcal{R}\lesssim (1+E)^{\frac{1-k}{2}}\lesssim 1,$$
as long as $k\geq 1$. In other words, in the most likely scenario, the weights offset the growth of the potential and  the operator behaves like a Maxwell molecule. For this reason, we provide an exhaustive treatment of Maxwell molecules, which is of independent interest. In fact, we improve and extend existing results for such operators. As mentioned, convolution estimates without additional angular cut-offs for Maxwell molecules have been proved in \cite{AlCa10,AlCaGa10}. However, in these works the implicit constants in the estimates for $b=1$ are finite only when $r<3$. We opt for a self-contained, purely kinetic approach that allows us to extend the existing estimates on Maxwell molecules for $b$ satisfying \eqref{assumption 1}--\eqref{assumption 3} to arbitrarily large $r$ (even $r=\infty$), as long as $p<3$. This is the content of Theorem \ref{Maxwell theorem general}. We stress that our technique is also simpler, and only relies on elementary kinetic tools. Moreover, most of the aforementioned papers that introduce angular cut-offs in the gain operator fall into this category, since the cut-offs allow only energy-distributing collisions to occur. Our method could be used to revisit and extend these works as well.

\subsubsection*{Energy-absorbing collisions} 
\begin{figure}[h!]
    \centering
    \includegraphics[width=0.4\textwidth]{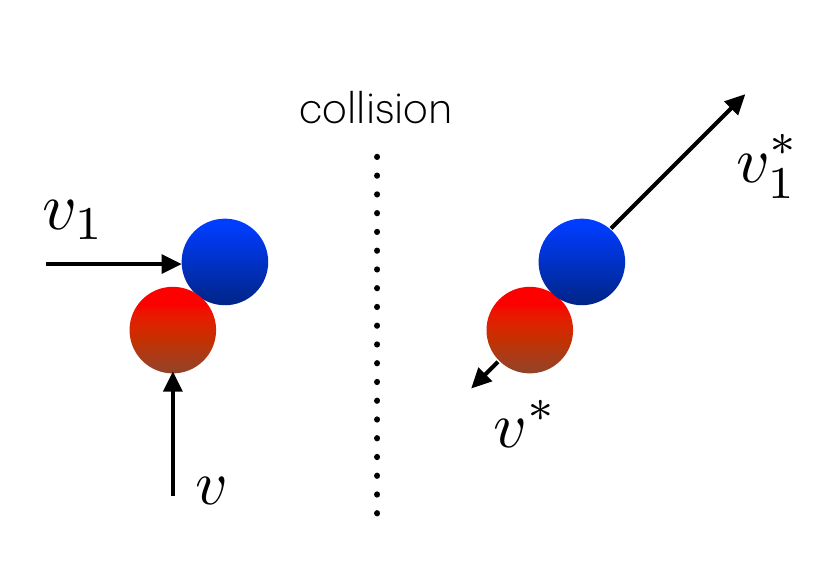} 
    \caption{Energy-absorbing collision with $|v_1|\gtrsim |v|$,  $v^*_1\approx v+v_1$}
    \label{energy absorbing figure 1}
\end{figure}

\begin{figure}
\centering
    \includegraphics[width=0.4\textwidth]{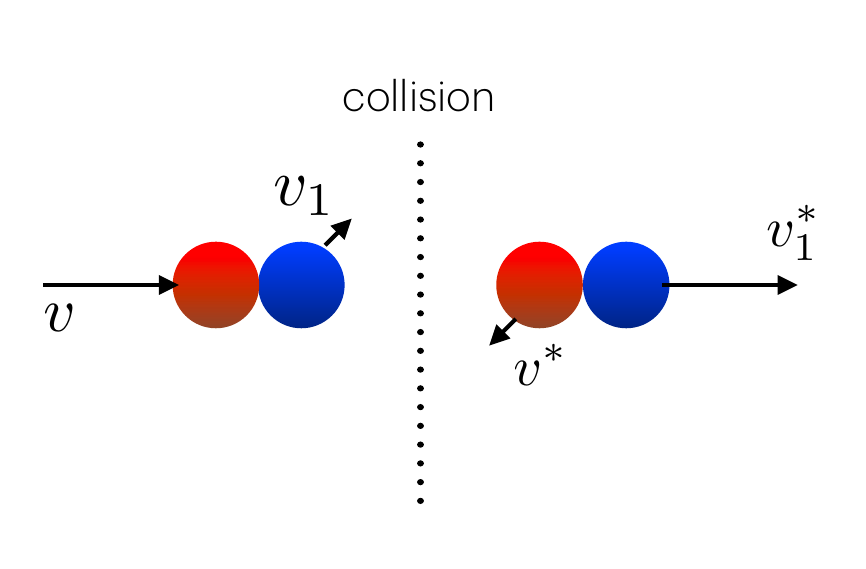} 
    \caption{Energy-absorbing collisions with $|v_1|\ll|v|$, $v_1^*\approx v$}
    \label{energy absorbing figure 2}
\end{figure}

The second, more pathological scenario consists of one of the outgoing particles accumulating almost all the energy of the collision i.e. $|v_1^*|\approx E^{1/2}$, $v^*\approx 0$. We call those {\it energy-absorbing collisions}. They distinguish hard potentials from Maxwell molecules. Following the billiard analogy, they can be visualized as two balls colliding, one staying still after the collision, while the other absorbs all of the energy. In Figures \ref{energy absorbing figure 1}, \ref{energy absorbing figure 2}, we provide examples of the two prototype energy-absorbing  collisions, depending on whether $|v_1|\gtrsim |v|$ or $|v_1|\ll |v|$.


It is easy to see that energy-absorbing collisions happen with non-trivial probability. For example when $\sigma\approx \frac{2V}{|u|}$, the collisional law \eqref{post-collisional velocities} yields 
$v^*\approx 0,\quad v_1^*\approx V$. However, in this case $\mathcal{R}\approx |u|=|v^*-v_1^*|\approx |v_1^*|\approx E^{1/2}$, where for the second equality we used \eqref{conservation of relative velocity}. This leads to a growth of moments which cannot be controlled point-wise. Hence a crude bound would result in a moment-increasing estimate, which is not sufficient to understand the time evolution of the system, even locally in time. 

Therefore to understand the mapping properties of the gain operator, it is necessary to  identify precisely where these pathological behaviors happen in $(v_1,\sigma)$ space. Mathematically, they correspond to $(v_1,\sigma)$ pairs such that $\mathcal{R}$ is not bounded as $\vert v \vert \to +\infty.$ The main contribution of this paper is to identify and control the behavior of the operator in these regions. It is captured by lower bounds on the small post-collisional velocity, see Lemma \ref{geometric}. The singular directions correspond to $(v_1,\sigma)$ configurations where said lower bound degenerates. The second challenge is to estimate the overall contribution of these regions for the gain operator. This presents an additional challenge given that the lower bounds obtained above are complicated nonlinear functions of $v_1$ and $\sigma.$ To overcome this difficulty, we use Bobylev variables (see Section \ref{sec:toolbox}). However, with this change of coordinate being singular, we cannot uniformly cover the entire set of $(v_1,\sigma)$ pairs. Therefore, we distinguish two cases based on whether $v_1^*$ and $v-v_1^*$ are closer to aligned or perpendicular. In the closer to aligned case we use the lower bound \eqref{lower bound for nonperp}, while in the closer to perpendicular case we use the lower bound \eqref{lower bound for perp}.
\color{black}

The gain operator contains a crucial angular averaging component: rather than the raw magnitude of $\mathcal{R}$, it is its expectation with respect to the angular variable that is relevant for the mapping properties. Therefore, the key mathematical quantities to consider are the averages of $\mathcal{R}$ (or more precisely its lower bounds derived in Lemma \ref{geometric}) over $\mathbb{S}^2$ (or in regions near singular directions). These crucial spherical integral are performed in Lemmas \ref{lem: lemma on A} and \ref{lem: perp lemma} for the aligned and perpendicular scenarios respectively.

Ultimately, putting all these ingredients together, we show that the gain operator satisfies the same estimates as the Maxwell molecule operator, with suitable weights, see Theorem \ref{main theorem}.

\subsection{Main results}

We now state the main result of this paper:

\begin{theorem}\label{main theorem} Let $0<\gamma <2 $ and $b$ satisfying \eqref{assumption 1}--\eqref{assumption 3}. Let $1\leq p,q<3$ with $p\geq \frac{2}{2-\gamma}$, and $\frac{2}{2-\gamma}\leq r\leq \infty$, such that 
$\frac{1}{p}+\frac{1}{q}=1+\frac{1}{r}$, and fix $k\geq l>\frac{2}{p'}$ with $k>\frac{1}{q'}$. Then, $Q^+$ extends uniquely to a  bounded, bilinear  operator $Q^+:\la v\ra^{-k}L^p_v\times\la v\ra^{-k}L^q_v\to \la v\ra^{-k}L^r_v$,
which satisfies the moment-preserving estimate
\begin{equation}\label{main estimate}
\|\l v\r^{k} Q^+(f,g)\|_{L^r}\lesssim  \|\l v\r^l f\|_{L^{p}}\|\l v\r^k g\|_{L^{q}}. 
\end{equation}
\end{theorem}

\begin{remark}
As it will become clear from the proof, the exponents range and the weights threshold are most likely optimal (although we do not prove it), at least within the framework of our approach.    
\end{remark}

\begin{remark}\label{remark on exponential}
Unlike for polynomial weights, we do not expect moment-preserving bounds on $Q^+$ to hold. Indeed, consider Gaussian weights of the form $e^{|v|^2}$. To prove the corresponding weighted estimates, one would need to control the ratio
$$\mathcal{R}_{exp}:=\frac{|u|^\gamma e^{|v|^2}}{e^{|v^*|^2}e^{|v_1^*|^2}}=|u|^\gamma e^{|v|^2-|v^*|^2-|v_1^*|^2}=|u|^\gamma e^{-|v_1|^2},$$
by conservation of energy. Now for $|v|\gg1$  and $|v_1|\ll1$, we have $\mathcal{R}_{exp}\approx |v|^\gamma$. However, unlike for polynomial weights, this quantity does not depend on the scattering angle $\sigma.$ Therefore collisional averaging cannot offset the growth induced by the potential.
\end{remark}

\subsection*{Organization of the paper} The paper is organized as follows:
\begin{itemize}[]
    \item[--] In Section  \ref{sec:toolbox} \color{black} we collect the main technical tools needed in the analysis, including Bobylev variables, custom collisional averaging lemmas as well as computations of two spherical integrals.
    
    \item[--] In Section \ref{sec:maxwell}, we prove the convolution estimates for Maxwell molecules, namely Theorem \ref{Maxwell theorem general}. This result in itself improves on the existing literature on the subject. It is also needed for non-singular part of the hard-sphere gain operator in the following section. \color{black}
    \item[--] In Section \ref{sec:hard}, we prove Theorem \ref{main theorem} which provides moment-preserving convolution estimates for hard potentials, and is the heart of our contribution. For the singular part,  we crucially rely on tools from Section \ref{sec:toolbox} and further develop them,  while the non-singular part follows from the  Maxwell molecules estimate from Section \ref{sec:maxwell}. \color{black}
\end{itemize}
\subsection*{Notation} Throughout this paper, given $v,v_1\in \R^3$ and $\sigma\in \S^2$, we will write 
\begin{align}\label{notation summary}
    v^* = V - \frac{\vert u \vert}{2} \sigma,  \ \ \ \ v_1^*=V+\frac{|u|}{2}\sigma,\ \ \ \  V = \frac{v+v_1}{2} , \ \ \ u = v-v_1, \ \ \ E=|v|^2+|v_1|^2 .
\end{align}

We will also use the following notation:

\begin{itemize}
\item We write $C_c(\R^3)$ for the space of continuous, compactly supported functions in $\R^3$,  $C_c^\infty(\R^3)$ for the space of smooth, compactly supported functions and $C_b(\R^3)$ for the space of continuous and bounded functions. 
\item We write $A \lesssim B,$ to mean that there exists a numerical constant $C > 0$, independent of $A,B$, such that $A \leq C B.$ 
\item We write $A \approx B$ to signify that both $A \lesssim B$ and $B \lesssim A$ hold.
\item We use the standard Japanese bracket notation: $\langle v \rangle : = \sqrt{1 + \vert v \vert^2},$ where $\vert \cdot \vert$ denotes the $\ell_v^2$--norm of the vector $v \in \mathbb{R}^3.$ 
\item For $v \in \mathbb{R}^3$ with $v\neq 0$,  we denote $\widehat{v} : = v/\vert v \vert.$
\end{itemize}

\medskip
\noindent{\bf Acknowledgments}  The authors are  grateful to the Institute for Advanced Study for providing a peaceful environment to conduct this research during July 2025. I.A.  was supported by the NSF grant No. DMS-2418020 and the PSC-CUNY Research Award 68653-0056.




\section{Technical toolbox}\label{sec:toolbox}
 In this section, we present some technical results fundamental to our analysis. First, we introduce the so-called Bobylev variables \cite{Bo75, Bo88}, and outline the properties of this coordinate system in Lemma \ref{chgvarkin}. Then, we  recall the  classical involutionary pre-post collisional change of variables in Lemma \ref{lem:involution}. We finally conjuct all the above by proving a custom useful  singular integral estimate in  Lemma \ref{lem:singular involutionary estimate}.
\color{black}
\subsection{Bobylev variables}
Given $\sigma\in\S^2$, we define the maps $R_{\sigma}^+,R_{\sigma}^-:\R^3\to\R^3$  by
\begin{align}
R_\sigma^+(y)&=\frac{y}{2}+\frac{|y|}{2}\sigma,\label{R+ def}\\
R_\sigma^-(y)&=\frac{y}{2}-\frac{|y|}{2}\sigma.\label{R- def}
\end{align}
One can readily verify that the following relations hold for any $y\in\R^3$:
\begin{align}
R_\sigma^+(y)+R_\sigma^-(y)&=y,\label{decomposition of u}\\
R_\sigma^+(y)\cdot R_\sigma^-(y)&=0,\label{orthogonality}\\
|R_\sigma^+(y)|^2+|R_\sigma^-(y)|^2&=|y|^2.\label{pythagorean}
\end{align} 

Note that given $v,v_1\in\R^3$, $\sigma\in\S^2$, we can connect $R_\sigma^+, R_\sigma^-$ with $v^*,v_1^*$ as follows:
\begin{align}
v^*&=v-R_\sigma^+(u),\label{R+ v*}\\
v_1^*&=v-R_\sigma^-(u).\label{R- v_1*}
\end{align}

We now state some finer substitution properties of $R_\sigma^+, R_\sigma^-$ that will be important in quantifying the collisional averaging behavior of the gain operator. The proof of the following result can be found in e.g. the work of Arsenio \cite{Arsenio}, as well as our previous work \cite{AmLe24}. For convenience of the reader we present a detailed proof.

\begin{proposition} \label{chgvarkin}
Let $\sigma\in\S^{2}$ and $\epsilon \in\{+,-\}$.
Then the map
$$R_\sigma^{\epsilon}:\lbrace y\in\R^3: y \cdot \sigma \neq - \epsilon  \vert y\vert \rbrace\to \lbrace \nu\in\R^3: \epsilon \,(\nu \cdot \sigma) > 0 \rbrace,$$ is a diffeomorphism with inverse
\begin{align}\label{inverse function}
    \big(R_{\sigma}^{\epsilon}\big)^{-1}(\nu) = 2 \nu  -\frac{\vert \nu \vert}{ (\widehat{\nu}\cdot \sigma )} \sigma,
\end{align}
and Jacobian 
\begin{equation}\label{Jacobian}
    \jac\,(R_\sigma^{\epsilon})^{-1}(\nu)=\frac{4 }{(\widehat{\nu}\cdot\sigma)^2}.
\end{equation}    
Moreover, for any $y\in\R^3$ with $y\cdot\sigma\neq -\epsilon|y|$, we have
\begin{align}
|R_\sigma^\epsilon(y)\cdot\sigma|&=\epsilon\,(R_\sigma^\epsilon(y)\cdot\sigma),\label{sign}\\
|y|&=\frac{|R_{\sigma}^{\epsilon}(y)|}{|\widehat{R}_{\sigma}^{\epsilon}(y)\cdot\sigma|}\label{magnitude},\\
\widehat{y} \cdot \sigma &= \epsilon \left(2 |\widehat{R}_{\sigma}^{\epsilon}(y)\cdot \sigma|^2 -1\right). \label{angle}
\end{align}
Finally for $y\in \R^3$ with $y\cdot\sigma\neq\pm |y|$, we have
\begin{equation}\label{R+ R- relation}
 |\widehat{R}_\sigma^+(y)|^2+|\widehat{R}_\sigma^-(y)|^2=1. 
\end{equation}
\end{proposition}
\begin{proof}
Let $\sigma\in \S^2$ and $\epsilon\in\{+,-\}$. Let us also denote
\begin{align*}
A_\sigma^\epsilon&=\{y\in\R^3,\:\, y\cdot\sigma\neq -\epsilon|y|\},\\
B_\sigma^\epsilon&=\{\nu\in\R^3\,:\, \epsilon\, (\nu\cdot\sigma)>0\}.
\end{align*}

First, for $y\in A_{\sigma}^\epsilon$ we clearly have $R_\sigma^\epsilon(y)\cdot\sigma\neq 0$. Moreover, since $\epsilon\in\{+,-\}$, we obtain 
\begin{equation}\label{well defined map}
\begin{split}\epsilon(R_\sigma^\epsilon(y)\cdot\sigma)&=\epsilon\left(\frac{y\cdot\sigma}{2}+\epsilon\frac{|y|}{2}\right)=\frac{|y|}{2}(1+\epsilon\,\widehat{y}\cdot\sigma)\\
&=\frac{|y|}{2}|1+\epsilon\,\widehat{y}\cdot\sigma|=\frac{|y|}{2}|\epsilon+\widehat{y}\cdot\sigma|\\
&=\left|\frac{y\cdot\sigma}{2}+\epsilon\frac{|y|}{2}\right|=|R_\sigma^\epsilon(y)\cdot\sigma|,
\end{split}
\end{equation}
thus $R_\sigma^{\epsilon}: A_{\sigma}^\epsilon\to B_\sigma^\epsilon$,
and moreover \eqref{sign} holds.

In order to show that $R_{\sigma}^\epsilon$ is a diffeomorphism with inverse given by \eqref{inverse function}, consider $y\in A_{\sigma}^\epsilon$, $\nu\in B_\sigma^\epsilon$  with $\nu=R_\sigma^\epsilon(y)$. Then $y=2\nu-\epsilon|y|\sigma$, thus
$$|y|^2=4|\nu|^2-4\epsilon|y|(\nu\cdot\sigma)+|y|^2\Rightarrow |y|=\frac{|\nu|}{\epsilon(\widehat{\nu}\cdot\sigma)}, $$
and \eqref{inverse function} follows.   Moreover, it is clear that $R_\sigma^\epsilon,(R_\sigma^\epsilon)^{-1}$ are continuously differentiable in $A_\sigma^\epsilon, B_\sigma^\epsilon$ respectively, thus    
 \eqref{inverse function} follows.   Moreover, it is clear that $R_\sigma^\epsilon,\,(R_\sigma^\epsilon)^{-1}$ are continuously differentiable in $A_\sigma^\epsilon, B_\sigma^\epsilon$, thus $R_\sigma^\epsilon:A_\sigma^\epsilon\to B_\sigma^\epsilon$ is a diffeomorphism.

To compute the Jacobian, we differentiate \eqref{inverse function} in $\nu$, to obtain
$$D (R_{\sigma})^{-1}(\nu)=2 I_3-\sigma\nabla^Tf(\nu), $$
where $I_3$ is the identity matrix and $f(\nu)=\frac{|\nu|^2}{(\nu\cdot\sigma)}$. Thus, we have\footnote{we use the Linear Algebra identity $\det(\lambda I_n+vw^T)=\lambda^n(1+\lambda^{-1}v\cdot w)$, where $v,w\in\R^n$, $n\in\mathbb{N}$ and $\lambda\neq 0$.}
\begin{equation}\label{Jacobian proof}
    \jac (R_\sigma)^{-1}(\nu)=\det\left(2 I_3-\sigma\nabla^T f(\nu)\right)=8\left(1-\frac{1}{2}\nabla f(\nu)\cdot\sigma\right).
\end{equation}
We readily compute
$$\nabla f(\nu)\cdot\sigma=\frac{2(\nu\cdot\sigma)^2-|\nu|^2}{(\nu\cdot\sigma)^2}=2-\frac{|\nu|^2}{(\nu\cdot\sigma)^2},$$
which, combined with \eqref{Jacobian proof}, gives \eqref{Jacobian}.

It remains to prove \eqref{magnitude}--\eqref{R+ R- relation}.
Indeed, since $y=2R_\sigma^\epsilon(y)-\epsilon|y|\sigma$, we have 
$$|y|^2=\Big|2R_\sigma^\epsilon(y)-\epsilon|y|\sigma\Big|^2=4|R_\sigma^\epsilon(y)|^2+|y|^2-4\epsilon|y| (R_\sigma^\epsilon(y)\cdot\sigma),$$
thus using \eqref{sign} we obtain
$$|R_\sigma^\epsilon(y)|^2=\epsilon |y| (R_\sigma^\epsilon(y)\cdot\sigma)=|y|\,|R_\sigma^\epsilon(y)\cdot\sigma|=|y|\,|R_\sigma^\epsilon(y)|\,|\widehat{R}_\sigma^\epsilon(y)\cdot\sigma|,$$
and \eqref{magnitude} follows.

Now using \eqref{magnitude}, \eqref{sign},  we obtain
\begin{align*}
y\cdot\sigma&=2(R_\sigma^\epsilon(y)\cdot\sigma)-\epsilon |y|=\epsilon|y|\left(2\frac{\epsilon(R_\sigma^\epsilon(y)\cdot\sigma)}{|y|}-1\right)\\
&=\epsilon|u|\left(2\frac{|R_\sigma^\epsilon(y)\cdot\sigma|}{|y|}-1\right)=\epsilon|y|\left(2|\widehat{R}_\sigma^\epsilon(y)\cdot\sigma|^2-1\right),
\end{align*}
and \eqref{angle} follows.

Finally, for $y\in\R^3$ with $y\cdot\sigma\neq \pm |y|$ we obtain \eqref{R+ R- relation} by combining \eqref{pythagorean} with \eqref{magnitude}. 
\end{proof}

\subsection{Pre-post collisional change of variables}
Finally, we recall the well-known pre-post collisional change of variables, see e.g. \cite{Villani-review}, and record a singular estimate which will be useful in the $L^r$-- estimates. We provide the proof for convenience of the reader.
\begin{lemma} \label{lem:involution}
 The map $T: (v,v_1,\sigma)\in\R^3\times\R^3\times\S^2\to (v^*,v_1^*,\eta)\in\R^3\times\R^3\times\S^2$ where
 \begin{equation}
 \begin{cases}
 v^*=\frac{v+v_1}{2}-\frac{|v-v_1|}{2}\sigma\\
 v_1^*=\frac{v+v_1}{2}+\frac{|v-v_1|}{2}\sigma\\
 \eta=-\widehat{u}=\frac{v_1-v}{|v-v_1|}
 \end{cases}
\end{equation}
 is an involution of $\R^3\times\R^3\times\S^2$. Moreover, for any non-negative and continuously differentiable function $F:\R^3\times\R^3\times [-1,1]\to \R_+$,  there holds the change of variables formula
 \begin{equation}\label{change of variables formula}
\int_{\R^6\times\S^2}F(v^*,v_1^*,\widehat{u}\cdot\sigma)\,d\sigma\,dv_1\,dv=\int_{\R^6\times\S^2} F(v,v_1,\widehat{u}\cdot\sigma)\,d\sigma\,dv_1\,dv.  
 \end{equation}
 \end{lemma}
 \begin{remark}
It is customary in literature to abbreviate notation and write $(v,v_1)\to (v^*,v_1^*)$ to indicate the above change of variables, however we stress that it is not correct to perform this change of variables for a fixed $\sigma\in\S^2$.     
\end{remark}
 \begin{proof}
Writing $T^2:=T\circ T$, we compute
     \begin{align*}
      T^2(v,v_1,\sigma)&=T(v^*,v_1^*,\eta)=\begin{pmatrix}
        \frac{v^*+v_1^*}{2}-\frac{|v^*-v_1^*|}{2}\eta\\
        \frac{v^*+v_1^*}{2}+\frac{|v^*-v_1^*|}{2}\eta\\
        \frac{v_1^*-v^*}{|v^*-v_1^*|}
      \end{pmatrix} =\begin{pmatrix}
        \frac{v+v_1}{2}-\frac{|v-v_1|}{2}\eta\\
        \frac{v+v_1}{2}+\frac{|v-v_1|}{2}\eta\\
        \frac{v_1^*-v^*}{|v-v_1|}
        \end{pmatrix}
        =\begin{pmatrix}
            v\\
            v_1\\
            \sigma
        \end{pmatrix} .
     \end{align*}
     Therefore $T^2 = Id$, so $T$ is an involution.
     
Second, writing $u^*:=v^*-v_1^*=-|u|\sigma$, we compute
\begin{align*}
   u^*\cdot\eta= -|u|\sigma\cdot(-\widehat{u})=u\cdot\sigma,
\end{align*}
thus by \eqref{conservation of relative velocity} we obtain $\widehat{u} \cdot \sigma = \widehat{u^*} \cdot \eta.$ Therefore
\begin{align*}
\int_{\R^6\times\S^2}F(v^*,v_1^*,\widehat{u}\cdot\sigma)\,d\sigma\,dv_1\,dv&= \int_{\R^6\times\S^2}F(v^*,v_1^*,\widehat{u^*}\cdot\eta)\,d\sigma\,dv_1\,dv\\
&=   \int_{\R^6\times\S^2}F(v,v_1,\widehat{u}\cdot\sigma)\,d\sigma\,dv_1\,dv,
\end{align*}
since $T$ is an involution and thus has unitary Jacobian.
\end{proof}

\begin{lemma}\label{lem:singular involutionary estimate}
Let $\beta<2$ and $f,g\in C_c(\R^3)$. Then the following estimate holds
\begin{equation}\label{singular involutionary estimate}
\int_{\R^3\times\R^3\times\S^2}\frac{|f(v^*)g(v_1^*)|}{|\hat{R}_\sigma^\epsilon(u)\cdot\sigma|^\beta}\,d\sigma\,dv_1\,dv\lesssim \|f\|_{L^1}\|g\|_{L^1},\quad \epsilon\in\{-1,1\}.    
\end{equation}
\begin{proof}
Using \eqref{angle}, we have $|\hat{R}_\sigma^\epsilon(u)\cdot\sigma|^2=\frac{1+\epsilon(\hat{u}\cdot\sigma)}{2}$, thus Lemma \ref{lem:involution} implies
\begin{align*}
 \int_{\R^3\times\R^3\times\S^2}\frac{|f(v^*)g(v_1^*)|}{|\hat{R}_\sigma^\epsilon(u)\cdot\sigma|^\beta}\,d\sigma\,dv_1\,dv&= \int_{\R^3\times\R^3\times\S^2}\frac{|f(v^*)g(v_1^*)|}{\big(1+\epsilon( \hat{u}\cdot\sigma)\big)^\beta/2}\,d\sigma\,dv_1\,dv \\
 &=\int_{\R^3\times\R^3}|f(v)g(v_1)|\int_{\S^2}\frac{1}{\big(1+\epsilon(\hat{u}\cdot\sigma)\big)^{\beta/2}}\,d\sigma\,dv_1\,dv
\end{align*}
Now, integrating in spherical coordinates, we obtain
$$\int_{\S^2}\frac{1}{\big(1+\epsilon(\hat{u}\cdot\sigma)\big)^{\beta/2}}\,d\sigma\approx\int_{-1}^1 \frac{1}{(1+\epsilon x)^{\beta/2}}\,dx\lesssim 1,$$
since $\epsilon\in\{-1,1\}$ and $\beta<2$.

Therefore
$$\int_{\R^3\times\R^3\times\S^2}\frac{|f(v^*)g(v_1^*)|}{|\hat{R}_\sigma^\epsilon(u)\cdot\sigma|^\beta}\,d\sigma\,dv_1\,dv\lesssim \int_{\R^3\times\R^3}|f(v)g(v_1)|\,dv_1\,dv=\|f\|_{L^1}\|g\|_{L^1}.$$
\end{proof}

\end{lemma}

\section{Maxwell molecules} \label{sec:maxwell}
In this section, we prove Young inequalities for Maxwell molecules. The mapping properties in this case are better than for hard potentials, thus they provide a natural best-case scenario for this harder setting. Moreover, as we stressed in the introduction, most collisions in the hard potential case are energy distributing. Handling these collisions is equivalent to dealing with Maxwell molecules. 

In Section \ref{maxwell: averaging estimates} we prove some useful estimates of singular integrals involving Bobylev variables. They are then used in Sections \ref{Maxwell: Linfty} and \ref{Maxwell: Lr} to prove the desired mapping properties. 

\subsection{Weightless averaging estimates} \label{maxwell: averaging estimates}
Here, we use Lemma \ref{chgvarkin} to obtain the main averaging estimates when moments do not increase. 
\begin{lemma}\label{basic averaging lemma}
Let $\beta_1<1$, $\beta_2>1/2$ and $1\leq p<\infty$. For $f\in C_c(\R^3)$, we define
\begin{align}
 \mathcal{I}^1_{\beta_1,p}[f](v)&:=\int_{\R^3\times\S^2}\frac{|f(v^*)|^p\mathds{1}_{\widehat{u}\cdot\sigma>0}}{|\widehat{R}_\sigma^-(u)\cdot\sigma|^{2\beta_1}}\,d\sigma\,dv_1,\label{definition of I_1 Maxwell}\\
 \mathcal{I}^2_{\beta_2,p}[f](v)&:= \int_{\R^3\times\S^2} |f(v_1^*)|^p|\widehat{R}_\sigma^\pm(u)\cdot\sigma|^{2\beta_2}\,d\sigma\,dv_1\label{definition of I_2 Maxwell}.
\end{align}
Then, there hold the estimates
\begin{align}
\|\mathcal{I}_{\beta_1,p}^1[f]\|_{L^\infty_v}&\lesssim \|f\|_{L^p_v}^p,\label{averaging estimate on I_1 Maxwell}\\
\|\mathcal{I}_{\beta_2,p}^2[f]\|_{L^\infty_v}&\lesssim \|f\|_{L^p_v}^p.\label{averaging estimate on I_2 Maxwell}
\end{align}
\end{lemma}
\begin{proof} Consider $f\in C_c(\R^3)$. We first prove \eqref{averaging estimate on I_1 Maxwell}. Let us write $\chi(x):=\mathds{1}_{x>0}$
and fix $v\in\R^3$. Using \eqref{R+ v*}, \eqref{angle}, \eqref{R+ R- relation} and the substitution $y:=u=v-v_1$ we write
\begin{align*}
\mathcal{I}_{\beta_1,p}^1[f](v)&=\int_{\S^2}\int_{\R^3}\frac{|f(v-R_\sigma^+(u))|^p\chi(2|\widehat{R}_\sigma^+(u)\cdot\sigma|^2-1)}{\left(1-|\widehat{R}_\sigma^+(u)\cdot\sigma|^2\right)^{\beta_1}}\,dv_1\,d\sigma\\
&=\int_{\S^2}\int_{\R^3}\frac{|f(v-R_\sigma^+(y))|^p\chi(2|\widehat{R}_\sigma^+(y)\cdot\sigma|^2-1)}{\left(1-|\widehat{R}_\sigma^+(y)\cdot\sigma|^2\right)^{\beta_1}}\,dy\,d\sigma.
\end{align*}

Now, we use Proposition \ref{chgvarkin} to substitute $\nu:=R_\sigma^+(y)$ and  obtain
\begin{align*}
\mathcal{I}_{\beta_1,p}^1[f](v)&\approx\int_{\S^2}\int_{\R^3}\frac{|f(v-\nu)|^p\chi(2|\widehat{\nu}\cdot\sigma|^2-1)}{|\widehat{\nu}\cdot\sigma|^2\left(1-|\widehat{\nu}\cdot\sigma|^2\right)^{\beta_1}}\,d\nu\,d\sigma\\
&= \int_{\R^3} |f(v-\nu)|^p\int_{\S^2}\frac{\chi(2|\widehat{\nu}\cdot\sigma|^2-1)}{|\widehat{\nu}\cdot\sigma|^2\left(1-|\widehat{\nu}\cdot\sigma|^2\right)^{\beta_1}}\,d\sigma\,d\nu\\
&\approx\int_{\R^3}|f(v-\nu)|^p\left(\int_{0}^1\frac{\chi(2x^2-1)}{x^2(1-x^2)^{\beta_1}}\,dx\right)\,d\nu\\
&\lesssim \int_{\R^3} |f(v-\nu)|^p\left(\int_{1/\sqrt{2}}^1\frac{1}{(1-x^2)^{\beta_1}}\,dx\right)\,d\nu\\
&\lesssim \int_{\R^3} |f(v-\nu)|^p\left(\int_{1/\sqrt{2}}^1\frac{1}{(1-x)^{\beta_1}}\,dx\right)\,d\nu\\
&\lesssim \int_{\R^3}|f(v-\nu)|^p\,d\nu=\|f\|_{L^p_v}^p,
\end{align*}
where we used  the fact that $\beta_1<1$ for the convergence of the integral in $x$ and translation invariance for the last equality. Estimate \eqref{averaging estimate on I_1 Maxwell} is proved.

We now prove \eqref{averaging estimate on I_2 Maxwell}. Fix again $v\in\R^3$. Using \eqref{R- v_1*} and \eqref{angle}, we write
\begin{align*}
\mathcal{I}_{\beta_2,p}^2[f](v)&=\int_{\S^2}\int_{\R^3}|f(v-R_\sigma^-(u))|^p|\widehat{R}_\sigma^-(u)\cdot\sigma|^{2\beta_2}\,dv_1\,d\sigma\\
&=\int_{\S^2}\int_{\R^3}|f(v-R_\sigma^-(y))|^p|\widehat{R}_\sigma^-(y)\cdot\sigma|^{2\beta_2}\,dy\,d\sigma.
\end{align*}

Using again Proposition \ref{chgvarkin}  to substitute $\nu:=R_\sigma^-(y)$,  we obtain
\begin{align*}
\mathcal{I}_{\beta_2,p}^2[f](v)&\approx \int_{\S^2}\int_{\R^3}|f(v-\nu)|^p|\widehat{\nu}\cdot\sigma|^{2(\beta_2-1)}\,d\nu\,d\sigma\\
&=\int_{\R^3}|f(v-\nu)|^p\int_{\S^2}|\widehat{\nu}\cdot\sigma|^{2(\beta_2-1)}\,d\sigma\,d\nu\\
&\approx \int_{\R^3}|f(v-\nu)|^p\left(\int_{0}^1 x^{2(\beta_2-1)}\,dx\right)\,d\nu\\
&\approx\int_{\R^3}|f(v-\nu)|^p\,d\nu=\|f\|_{L^p_v}^p,
\end{align*}
where we used the fact that $\beta_2>1/2$ for the convergence of the integral in $x$. Estimate \eqref{averaging estimate on I_2 Maxwell} is proved.
\end{proof}

\begin{lemma}
 Let $1\leq p,q<\infty$ and $f,g\in C_c(\R^3)$. Then the following estimate holds
\begin{align}
\sup_{v\in\R^3}\int_{\R^3\times\S^2}\frac{|\hat{R}_\sigma^+(u)\cdot\sigma|^\alpha}{|\hat{R}_\sigma^{-}(u)\cdot\sigma|^\beta}|f(v^*)|^p\,d\sigma\,dv_1&\lesssim \|f\|_{L^p}^p,\qquad \alpha>1,\quad \beta<2.\label{Maxwell averaging 1}
\end{align}
Moreover, we have the truncated estimate
\begin{align}
\sup_{v\in\R^3}\int_{\R^3\times\S^2}\frac{\mathds{1}_{\hat{u}\cdot\sigma>0}}{|\hat{R}_\sigma^{-}(u)\cdot\sigma|^\beta}|f(v^*)|^p\,d\sigma\,dv_1&\lesssim \|f\|_{L^p}^p,\qquad \beta<2, \label{Maxwell averaging truncated 1}.
\end{align}
\end{lemma}
\begin{proof}
We first prove \eqref{Maxwell averaging 1}. Let $\alpha>1$, $\beta<2$, and fix $v\in\R^3$. Write 
$$\mathcal{I}_{\alpha,\beta}:=\int_{\R^3\times\S^2}\frac{|\hat{R}_\sigma^+(u)\cdot\sigma|^\alpha}{|\hat{R}_\sigma^{-}(u)\cdot\sigma|^\beta}|f(v^*)|^p\,d\sigma\,dv_1.$$
Using  \eqref{R+ v*}, \eqref{R+ R- relation}, and the substitution $y:=u=v-v_1$, we have 
\begin{align*}
\mathcal{I}_{\alpha,\beta}&=\int_{\R^3\times\S^2} \frac{|\hat{R}_\sigma^+(u)\cdot\sigma|^\alpha}{\left(1-|\hat{R}_\sigma^+(u)\cdot\sigma|^2\right)^{\beta/2}}|f(v-R_\sigma^+(u))\,d\sigma\,dv_1\\
&=\int_{\S^2}\int_{\R^3} \frac{|\hat{R}_\sigma^+(y)\cdot\sigma|^\alpha}{\left(1-|\hat{R}_\sigma^+(y)\cdot\sigma|^2\right)^{\beta/2}}|f(v-R_\sigma^+(y))\,dy\,d\sigma.
\end{align*}
For fixed $\sigma\in\S^2$, we apply Proposition \ref{chgvarkin} to substitute $\nu:=R_\sigma^+(y)$ and obtain
\begin{align*}
\mathcal{I}_{\alpha,\beta}&=\int_{\S^2}\int_{\R^3}\frac{|f(v-\nu)|^p}{|\hat{\nu}\cdot\sigma|^{2-\alpha}\left(1-|\hat{\nu}\cdot\sigma|^2\right)^{\beta/2}}\,d\nu\,d\sigma\\
&=\int_{\R^3}|f(v-\nu)|^p\int_{\S^2}\frac{1}{|\hat{\nu}\cdot\sigma|^{2-\alpha}\left(1-|\hat{\nu}\cdot\sigma|^2\right)^{\beta/2}}\,d\sigma\,d\nu.
\end{align*}
Integrating in spherical coordinates, we estimate
\begin{align*}
\int_{\S^2}& \frac{1}{|\hat{\nu}\cdot\sigma|^{2-\alpha}\left(1-|\hat{\nu}\cdot\sigma|^2\right)^{\beta/2}}\,d\sigma \approx\int_0^1 \frac{1}{x^{2-\alpha}(1-x^2)^{\beta/2}}\,dx\\
&\lesssim\int_0^{1/2}x^{\alpha-2}\,dx+\int_{1/2}^1\frac{1}{(1-x^2)^{\beta/2}}\,dx\\
&\leq \int_0^{1/2}x^{\alpha-2}\,dx+\int_{1/2}^1\frac{1}{(1-x)^{\beta/2}}\,dx\lesssim 1,
\end{align*}
since $\alpha>1$ and $\beta<2$. Therefore,
$$\mathcal{I}_{\alpha,\beta}\lesssim \int_{\R^3}|f(v-\nu)|^p\,d\nu=\|f\|_{L^p}^p,$$
by translation invariance. Estimate \eqref{Maxwell averaging 1} is proved.

We know prove \eqref{Maxwell averaging truncated 1}. Let us write $\chi(x):=\mathds{1}_{x>0}$
and 
$$\widetilde{\mathcal{{I}}}_\beta:=\int_{\R^3\times\S^2}\frac{\chi(\hat{u}\cdot\sigma)}{|\hat{R}_\sigma^-(u)\cdot\sigma|^\beta}|f(v^*)|^p\,d\sigma\,dv_1.$$
Using \eqref{R+ v*}, \eqref{angle}, \eqref{R+ R- relation} and the substitution $y:=u=v-v_1$ we obtain
\begin{align*}
\widetilde{\mathcal{I}}_\beta&=\int_{\S^2}\int_{\R^3}\frac{\chi(2|\hat{R}_\sigma^+(u)\cdot\sigma|^2-1)}{\left(1-|\hat{R}_\sigma^+(u)\cdot\sigma|^2\right)^{\beta/2}}|f(v-R_\sigma^+(u))|^p\,dv_1\,d\sigma\\
&=\int_{\S^2}\int_{\R^3}\frac{\chi(2|\hat{R}_\sigma^+(y)\cdot\sigma|^2-1)}{\left(1-|\hat{R}_\sigma^+(y)\cdot\sigma|^2\right)^{\beta/2}}|f(v-R_\sigma^+(y))|^p\,dy\,d\sigma.
\end{align*}
For fixed $\sigma\in\S^2$, we apply Proposition \ref{chgvarkin} to substitute $\nu:=R_\sigma^+(y)$ and obtain
\begin{align*}
\widetilde{\mathcal{I}}_\beta&=\int_{\S^2}\int_{\R^3}\frac{\chi(2|\hat{\nu}\cdot\sigma|^2-1)}{|\hat{\nu}\cdot\sigma|^{2}\left(1-|\hat{\nu}\cdot\sigma|^2\right)^{\beta/2}}|f(v-\nu)|^p\,d\nu\,d\sigma\\
&=\int_{\R^3}|f(v-\nu)|^p\int_{\S^2}\frac{\chi(2|\hat{\nu}\cdot\sigma|^2-1)}{|\hat{\nu}\cdot\sigma|^{2-\alpha}\left(1-|\hat{\nu}\cdot\sigma|^2\right)^{\beta/2}}\,d\sigma\,d\nu.
\end{align*}
Integrating in spherical coordinates, we compute
\begin{align*}\int_{\S^2}&\frac{\chi(2|\hat{\nu}\cdot\sigma|^2-1)}{|\hat{\nu}\cdot\sigma|^{2-\alpha}\left(1-|\hat{\nu}\cdot\sigma|^2\right)^{\beta/2}}\,d\sigma\approx \int_0^1 \frac{\chi(2x^2-1)}{x^2(1-x^2)^{\beta/2}}\,dx=\int_{|x|>1/\sqrt{2}}\frac{1}{x^2(1-x^2)^{\beta/2}}\,dx\\
&\approx \int_{|x|>1/\sqrt{2}}\frac{1}{(1-x^2)^{\beta/2}}\,dx\lesssim \int_{|x|>1/\sqrt{2}}\frac{1}{(1-x)^{\beta/2}}\,dx\lesssim 1,
\end{align*}
since $\beta<2$. Therefore,
$$\widetilde{\mathcal{I}}_{\beta}\lesssim \int_{\R^3}|f(v-\nu)|^p\,d\nu=\|f\|_{L^p}^p,$$
by translation invariance.

\end{proof}

\subsection{The Maxwell molecules operator}
In this section, we perform a simple symmetry-based restriction to a hemisphere in the gain operator for Maxwell molecules. We also state the main result of this Section.

Given $b$ satisfying \eqref{assumption 1}--\eqref{assumption 3}, we define the Maxwell bilinear collisional operator with angular cross-section $b$ as
\begin{equation}\label{Maxwell notation}
Q_{M,b}^+(f,g)=\int_{\R^3\times\S^2}b(\widehat{u}\cdot\sigma)f(v^*)g(v_1^*)\,d\sigma\,dv_1 .   
\end{equation}
Using spherical symmetry we can write  
\begin{align}
Q_M^+(f,g)(v)&=\int_{\R^3\times\S^2}b(\hat{u}\cdot\sigma)f(v^*)g(v_1^*)\mathds{1}_{\hat{u}\cdot\sigma>0}\,d\sigma\,dv_1+\int_{\R^3\times\S^2}b(\hat{u}\cdot\sigma)f(v^*)g(v_1^*)\mathds{1}_{\hat{u}\cdot\sigma<0}\,d\sigma\,dv_1\notag\\
&=\int_{\R^3\times\S^2}b(\hat{u}\cdot\sigma)f(v^*)g(v_1^*)\mathds{1}_{\hat{u}\cdot\sigma>0}\,d\sigma\,dv_1+\int_{\R^3\times\S^2}b(-\hat{u}\cdot\sigma)g(v^*)f(v_1^*)\mathds{1}_{\hat{u}\cdot\sigma>0}\,d\sigma\,dv_1\notag\\
&=\int_{\R^3\times\S^2}b(\hat{u}\cdot\sigma)\Big(f(v^*)g(v_1^*)+g(v^*)f(v_1^*)\Big)\mathds{1}_{\hat{u}\cdot\sigma>0}\,d\sigma\,dv_1\label{Maxwell truncated}
\end{align}
where for the second to last line we substituted $\sigma\to -\sigma$ in the second term of the summation (which exchanges $v^*$ and $v_1^*$ by \eqref{post-collisional velocities}), and to obtain the last line we used \eqref{assumption 2}.

As mentioned, when $b$ is constant, the previous works \cite{AlCa10,AlCaGa10} only provide an estimate when $r<3$. Here, we are able to extend the estimates to arbitrary $r$, as long as $p,q<3$. Namely, we prove the following:
\begin{theorem}\label{Maxwell theorem general} 
Let $b$ satisfying \eqref{assumption 1}--\eqref{assumption 3}. Let $1\leq p,q,r\leq\infty$ with $1\leq p,q<3$  such that $\frac{1}{p}+\frac{1}{q}=1+\frac{1}{r}$. 
Then, $Q_{M,b}^+$ extends uniquely  to a bounded, bilinear operator $Q_{M,b}^+:L^p\times L^q\to L^r$ satisfying the estimate
 \begin{equation}\label{Maxwell general estimate}
\|Q^+_{M,b}(f,g)\|_{L^r}\lesssim \|f\|_{L^p}\|g\|_{L^q}.     
 \end{equation}
\end{theorem}
The proof of this Theorem will occupy the rest of this section. We will treat the cases $r=\infty$ and $r<\infty$ separately.

\subsection{The $L^\infty$-- estimate} \label{Maxwell: Linfty}
First, we prove the convolution estimate corresponding to $r=\infty$.
\begin{proposition}\label{Maxwell prop Linf} 
Let $b$ satisfying \eqref{assumption 1}--\eqref{assumption 3}. Let $1< p,q<3$  such that $\frac{1}{p}+\frac{1}{q}=1$. 
Then, there holds the  estimate
 \begin{equation}\label{Maxwell Linf estimate}
\|Q^+_{M,b}(f,g)\|_{L^\infty}\lesssim \|f\|_{L^p}\|g\|_{L^q},\qquad\forall\, f,g\in C_c(\R^3).     
 \end{equation}
\end{proposition}

\begin{proof}Let $f,g\in C_c(\R^3)$.
Using \eqref{Maxwell truncated} and \eqref{assumption 3}, we obtain 
\begin{equation}\label{bound by constant kernel Maxwell}
|Q^+_{M,b}(f,g)(v)|\lesssim \int_{\R^3\times\S^2}|f(v^*)||g(v_1^*)|\mathds{1}_{\hat{u}\cdot\sigma>0}\,d\sigma\,dv_1+ \int_{\R^3\times\S^2}|g(v^*)||f(v_1^*)|\mathds{1}_{\hat{u}\cdot\sigma>0}\,d\sigma\,dv_1
\end{equation}

Since $p<3$ we have $p<2q$ so we can pick $\alpha_1$ with $1<\alpha_1<2q/p$. Using H\"older's inequality, we obtain
\begin{align}
\int_{\R^3\times\S^2}&|f(v^*)||g(v_1^*)|\mathds{1}_{\hat{u}\cdot\sigma>0}\,d\sigma\,dv_1=\int_{\R^3\times\S^2} \left(\frac{|f(v^*)|\mathds{1}_{\hat{u}\cdot\sigma>0}}{|\hat{R}_\sigma^-(u)\cdot\sigma|^{\alpha_1/q}}\right)\left(|\hat{R}_\sigma^-(u)\cdot\sigma|^{\alpha_1/q}|g(v_1^*)|\right)\,d\sigma\,dv_1\notag\\
 &\leq \left(\int_{\R^3\times\S^2}\frac{\mathds{1}_{\hat{u}\cdot\sigma>0}}{|\hat{R}_\sigma^-(u)\cdot\sigma|^{\alpha_1 p/q}}|f(v^*)|^p\,d\sigma\,dv_1\right)^{1/p}\left(\int_{\R^3\times\S^2}|\hat{R}_\sigma^-(u)\cdot\sigma|^{\alpha_1}|g(v_1^*)|^q\,d\sigma\,dv_1\right)^{1/q}\notag\\
 &\lesssim \|f\|_{L^p}\|g\|_{L^q},\label{positive estimate Maxwell 1}
\end{align}
where for the last line we used \eqref{Maxwell averaging 1}, \eqref{Maxwell averaging truncated 1} since $1<\alpha_1<2q/p$.

Now, since  $q<3$ as well, we have $q<2p$, so we can pick $\alpha_2$ with $1<\alpha_2<2p/q$. A symmetric argument then yields
\begin{equation}
\int_{\R^3\times\S^2}|g(v^*)||f(v_1^*)|\mathds{1}_{\hat{u}\cdot\sigma>0}\,d\sigma\,dv_1\lesssim \|g\|_{L^q}\|f\|_{L^p}=\|f\|_{L^p}\|g\|_{L^q}.\label{positive estimate Maxwell 2}   
\end{equation}

Combining \eqref{bound by constant kernel Maxwell}, \eqref{positive estimate Maxwell 1}, \eqref{positive estimate Maxwell 2}, the claim follows.
\end{proof}

\subsection{The $L^r$-- estimate} \label{Maxwell: Lr}
Now, we prove the convolution estimates for  $1\leq r<\infty$.
\begin{proposition}\label{Maxwell prop Lr} 
Let $b$ satisfying \eqref{assumption 1}--\eqref{assumption 3}.  Let $1\leq p,q,r<\infty$ with $p,q<3$ such that $\frac{1}{p}+\frac{1}{q}=1+\frac{1}{r}$. 
Then, there holds the estimate
\begin{equation}\label{Maxwell estimate Lr}
\|Q_{M,b}^+(f,g)\|_{L_v^r}\lesssim \|f\|_{L_v^p}\|g\|_{L_v^q},\qquad \forall\,f,g\in C_c(\R^3).  
\end{equation}
\end{proposition}
\begin{proof} Let $f,g\in C_c(\R^3)$. Since $p<3$, we have $p'>3/2$, thus the exponents' relation implies $\frac{1}{2p'}<\frac{1}{q'}+\frac{1}{r}$.
We choose $\alpha_1>1$ such that
\begin{equation}\label{alpha condition Maxwell}
1<\alpha_1<2p'\left(\frac{1}{q'}+\frac{1}{r}\right).    
\end{equation}
Then  $1-\frac{2p'}{\alpha_1 q'}<\frac{2p'}{\alpha_1 r}$. Choose $s_1\in\R$ with $1-\frac{2p'}{\alpha_1 q'}<s_1<\frac{2p'}{\alpha_1 r}$. Then \begin{equation}\label{s condition Maxwell}
\frac{(1-s_1)\alpha_1 q'}{p'}<2,\qquad\frac{s_1\alpha_1 r}{p'}<2.
\end{equation}

Using the exponents' relation and H\"older's inequality we bound
\begin{align*}
\int_{\R^3\times\S^2}|f(v^*)||g(v_1^*)|\mathds{1}_{\hat{u}\cdot\sigma>0}\,d\sigma\,dv_1&\leq \int_{\R^3\times\S^2}\left(\frac{|f(v^*)|^{p/q'}\mathds{1}_{\hat{u}\cdot\sigma>0}}{|\hat{R}_\sigma^-(u)\cdot\sigma|^{(1-s_1)\alpha_1/p'}}\right)\left(|\hat{R}_\sigma^-(u)\cdot\sigma|^{\alpha_1/p'}|g(v_1^*)|^{q/p'}\right)\\  
&\hspace{4cm}\times \left(\frac{|f(v^*)|^{p/r}|g(v_1^*)|^{q/r}}{|\hat{R}_\sigma^-(u)\cdot\sigma|^{s_1\alpha_1 /p'}}\right)\,d\sigma\,dv_1\notag\\
&\leq \mathcal{I}_1^{1/q'}\mathcal{I}_2^{1/p'}
 \left(\int_{\R^3\times\S^2}\frac{|f(v^*)|^p|g(v_1^*)|^q}{|\hat{R}_\sigma^-(u)\cdot\sigma|^{s_1\alpha_1 r/p'}}\,d\sigma\,dv_1\right)^{1/r},
\end{align*}
where 
\begin{align*}
\mathcal{I}_1&= \int_{\R^3\times\S^2}\frac{|f(v^*)|^p\mathds{1}_{\hat{u}\cdot\sigma>0}}{|\hat{R}^-_\sigma(u)\cdot\sigma|^{(1-s_1)\alpha_1 q'/p'}}\,d\sigma\,dv_1,\\
\mathcal{I}_2&=\int_{\R^3\times\S^2}|\hat{R}_\sigma^-(u)\cdot\sigma|^\alpha |g(v_1^*)|^q\,d\sigma\,dv_1,
\end{align*}
Now since $\alpha_1>1$ and $\frac{(1-s_1)\alpha_1 q'}{p'}<2$  \eqref{Maxwell averaging 1}, \eqref{Maxwell averaging truncated 1} imply
\begin{align*}
  I_1\lesssim \|f\|_{L^p}^p,\qquad I_2\lesssim \|g\|_{L^q}^q,  
\end{align*}
thus 
$$\int_{\R^3\times\S^2}|f(v^*)||g(v_1^*)\mathds{1}_{\hat{u}\cdot\sigma}\,d\sigma\,dv_1\lesssim \|f\|_{L^p}^{p/q'}\|g\|_{L^q}^{q/p'}\left(\int_{\R^3\times\S^2}\frac{|f(v^*)|^p|g(v_1^*)|^q}{|\hat{R}_\sigma^-(u)\cdot\sigma|^{s_1\alpha_1 r/p'}}\,d\sigma\,dv_1\right)^{1/r}.$$
Raising the above bound to the $r-$th power, integrating in $v$, and using the exponents' relation, we obtain
\begin{align*}
\left\|\int_{\R^3\times\S^2}|f(v^*)||g(v_1^*)|\mathds{1}_{\hat{u}\cdot\sigma>0}\,d\sigma\,dv_1\right\|_{L^r_v}^r&\lesssim \|f\|_{L^p}^{rp/q'}\|g\|_{L^q}^{rq/p'}\int_{\R^3\times\R^3\times\S^2}\frac{|f(v^*)|^p\|g(v_1^*)|^q}{|\hat{R}_\sigma^-(u)\cdot\sigma|^{s_1\alpha_1 r/p'}}\,d\sigma\,dv_1\,dv\\
&= \|f\|_{L^p}^{r-p}\|g\|_{L^q}^{r-q}\int_{\R^3\times\R^3\times\S^2}\frac{|f(v^*)|^p\|g(v_1^*)|^q}{|\hat{R}_\sigma^-(u)\cdot\sigma|^{s_1\alpha_1 r/p'}}\,d\sigma\,dv_1\,dv.
\end{align*}
Since $\frac{s_1\alpha_1 r}{p'}<2$, we use Lemma \ref{lem:singular involutionary estimate} to estimate the remainder term by $\|f\|_{L^p}^p\|g\|_{L^q}^q$, therefore 
\begin{equation}\label{Maxwel bound positive Lr 1}
\left\|\int_{\R^3\times\S^2}|f(v^*)||g(v_1^*)|\mathds{1}_{\hat{u}\cdot\sigma>0}\,d\sigma\,dv_1\right\|_{L^r_v}\lesssim
\|f\|_{L^p}\|g\|_{L^q}.
\end{equation}

Since $q<3$ as well, using the same argument, we can pick $\alpha_2>1$ and $s_2\in \R$ with 
$$ \frac{(1-s_2)\alpha_2 p'}{q'}<2,\qquad \frac{s_2\alpha_2 r}{q'}<2.$$
Then, a symmetric argument to above yields
\begin{equation}\label{Maxwel bound positive Lr 2}
\left\|\int_{\R^3\times\S^2}|g(v^*)||f(v_1^*)|\mathds{1}_{\hat{u}\cdot\sigma>0}\,d\sigma\,dv_1\right\|_{L^r_v}\lesssim
\|g\|_{L^q}\|f\|_{L^p}=\|f\|_{L^p}\|g\|_{L^q}.
\end{equation}
Combining \eqref{bound by constant kernel Maxwell}, \eqref{Maxwel bound positive Lr 1}, \eqref{Maxwel bound positive Lr 2}, the claim follows.
\end{proof}

\subsection{Proof of Theorem \ref{Maxwell theorem general}} By Propositions \ref{Maxwell prop Linf}, \ref{Maxwell prop Lr}, for any $1\leq p,q,r\leq \infty$ with $p,q<3$, such that $\frac{1}{p}+\frac{1}{q}=1+\frac{1}{r}$, we have
\begin{equation}\label{cc bound Maxwell}
\|Q^+_{M,b}(f,g)\|_{L^r}\lesssim \|f\|_{L^p}\|g\|_{L^q},\qquad\forall\, f,g\in C_c(\R^3).    
\end{equation}
The claim then follows by a standard density argument.

\section{Moment-preserving Young's inequality for hard potentials}\label{sec:hard}
In this final section, we prove our main result, Theorem \ref{main theorem}, for the bilinear operator $Q^+(f,g)$ given in \eqref{gain operator bilinear}. We first establish the necessary estimates for smooth, compactly supported functions  for $r=\infty$ and $2<r<\infty$ respectively.

As mentioned in the introduction, the main difficulty is to handle energy-distributing collisions. We start by proving two lower bounds on the small post-collisional velocity in Lemma \ref{geometric}. They capture the behavior of $\mathcal{R}$ (defined in \eqref{ratio intro}) in singular scattering directions. Given the angular averaging present in the gain operator, the relevant quantity for boundedness is the expected value of $\mathcal{R}$ over scattering directions. Thus crucial spherical integral calculations are performed in Lemmas \ref{lem: lemma on A} and \ref{lem: perp lemma}. These ingredients are then combined in the remainder of the section to prove the main result Theorem \ref{main theorem}.

\subsection{Comparison of post-collisional velocities}

\begin{lemma} \label{geometric}
Let $v,v_1\in\R^3$ and $\sigma\in\S^2$. Then, the following estimates hold a.e. in $\R^3\times\S^2$
\begin{align}
\la v^*\ra^2&\geq\la v_1^*\ra^2\left(1-\zeta_{v_1^*}|\hat{v_1^*}\cdot\sigma|\right),\qquad \zeta_{v_1^*}=\frac{|v_1^*|^2}{\la v_1^*\ra^2}\label{lower bound for nonperp}\\
\la v^*\ra^2&\geq \frac{\la v_1^*\ra^2 |\hat{R}_\sigma^\epsilon(u)\cdot\sigma|^2-2|v_1^*||R_\sigma^\epsilon(u)||\hat{v_1^*}\cdot\sigma||\hat{R}_\sigma^\epsilon(u)\cdot\sigma|+|R_\sigma^\epsilon(u)|^2}{|\hat{R}_\sigma^\epsilon(u)\cdot\sigma|^2},\quad \epsilon\in\{-1,1\}.\label{lower bound for perp}
\end{align}
\end{lemma}
\begin{proof}
Using \eqref{post-collisional velocities}, we obtain $v^*=v_1^*-|u|\sigma$, thus
\begin{align}\label{v* compared to v}
|v^*|^2\geq |v_1^*|^2+|u|^2-2 |v_1^*||u| |\hat{v_1^*}\cdot\sigma|.
\end{align}
To prove \eqref{lower bound for nonperp}, we use \eqref{v* compared to v} and the basic inequality $2ab\leq a^2+b^2$, we obtain
\begin{align}\label{lower bound on post-velocity magnitude}
|v^*|^2&=(|v_1^*|^2+|u|^2)\left(1-\frac{2|v_1^*||u|}{|v_1^*|^2+|u|^2}|\hat{v_1^*}\cdot\sigma|\right)\geq |v_1^*|^2(1-|\hat{v_1^*}\cdot\sigma|).    
\end{align}
Hence, we have 
\begin{align}
\la v^*\r^2=1+|v^*|^2\geq 1+|v_1^*|^2(1-|\hat{v_1^*}\cdot\sigma|)=\la v_1^*\ra^2-|v_1^*|^2|\hat{v_1^*}\cdot\sigma|=\la v_1^*\ra^2(1-\zeta_{v_1^*}|\hat{v_1^*}\cdot\sigma|),    
\end{align}
where $\zeta_{v_1^*}=\frac{|v_1^*|^2}{\la v_1^*\ra^2}$.

To prove estimate \eqref{lower bound for perp}, we use  \eqref{v* compared to v} and  \eqref{magnitude} to obtain
$$|v^*|^2\geq | v_1^*|^2+\frac{|R_\sigma^\epsilon(u)|^2}{|\hat{R}_\sigma^\epsilon(u)\cdot\sigma|^2}-2|v_1^*||R_\sigma^\epsilon(u)|\frac{|\hat{v_1^*}\cdot\sigma|}{|\hat{R}_\sigma^\epsilon(u)\cdot\sigma|},$$
thus 
\begin{align*}
\la v^*\ra^2&\geq \la v_1^*\ra^2+\frac{|R_\sigma^\epsilon(u)|^2}{|\hat{R}_\sigma^\epsilon(u)\cdot\sigma|^2}-2|v_1^*||R_\sigma^\epsilon(u)|\frac{|\hat{v_1^*}\cdot\sigma|}{|\hat{R}_\sigma^\epsilon(u)\cdot\sigma|}\\
&=\frac{\la v_1^*\ra^2 |\hat{R}_\sigma^\epsilon(u)\cdot\sigma|^2-2|v_1^*||R_\sigma^\epsilon(u)||\hat{v_1^*}\cdot\sigma||\hat{R}_\sigma^\epsilon(u)\cdot\sigma|+|R_\sigma^\epsilon(u)|^2}{|\hat{R}_\sigma^\epsilon(u)\cdot\sigma|^2}.
\end{align*}
\end{proof}

\subsection{Spherical averaging estimates}

\begin{lemma}\label{lem: lemma on A}
Let $\ell>1$.  For $w,\nu\in\R^3\setminus\{0\}$ and $0<\lambda<1$, we define 
 \begin{equation}\label{definition of A}
\mathcal{A}_{\ell}(w,\nu,\lambda):= \mathds{1}_{|\hat{w}\cdot\hat{\nu}|\geq \sqrt{2}/2}\int_{\S^2}\frac{\mathds{1}_{|\hat{w}\cdot\sigma|\geq\sqrt{3}/2}}{|\hat{\nu}\cdot\sigma|^2\left(1-\lambda|\hat{w}\cdot\sigma|\right)^{\ell}}\,d\sigma.     
 \end{equation}
 Then, there holds the estimate
 \begin{equation}
 \mathcal{A}_{\ell}(w,\nu,\lambda)\leq     \frac{4\pi}{\lambda (\ell-1)}(1-\lambda)^{1-\ell}
 \end{equation}
\end{lemma}

\begin{proof}  
 Consider $w,\nu\neq 0$ with $\sqrt{2}/2\leq|\hat{w}\cdot\hat{\nu}|<1$. We parametrize the sphere by choosing the  $z-$axis in the direction of $w$ and the $x-$axis along the projection of $\nu$ on  $w^\perp$ (which is non-trivial because $|\hat{w}\cdot\hat{\nu}|\neq 1$). Let $\theta\in[0,\pi]$ be  the polar angle and $\phi\in[0,2\pi)$ the azimuthal angle. Then any $\sigma\in\S^2$ can be uniquely written as $\sigma=(\cos\phi\sin\theta,\sin\phi\sin\theta,\cos\theta)$. In particular, we have $\hat{w}\cdot\sigma=\cos\theta$ and $\hat{\nu}=(\sin\psi,0,\cos\psi)$ where $\psi=\arccos(\hat{w}\cdot\hat{\nu})$. 
We estimate
\begin{align*}
|\hat{\nu}\cdot\sigma|\mathds{1}_{|\widehat{w}\cdot\sigma|\geq\sqrt{3}/2}\mathds{1}_{|\hat{w}\cdot\hat{\nu}|\geq \sqrt{2}/2}&=|\cos\phi\sin\theta\sin\psi+\cos\theta\cos\psi|\mathds{1}_{|\cos\theta|\geq\sqrt{3}/2}\mathds{1}_{|\cos\psi|\geq\sqrt{2}/2}\notag\\
&\geq \left(|\cos\theta\cos\psi|-|\sin\theta\sin\psi|\right)\mathds{1}_{|\cos\theta|\geq\sqrt{3}/2}\mathds{1}_{|\cos\psi|\geq\sqrt{2}/2}\\
&\geq \frac{\sqrt{2}(\sqrt{3}-1)}{4}>\frac{1}{4}.   
\end{align*}
Using the above bound, and integrating in spherical coordinates, we obtain
\begin{align*}
\mathcal{A}_\ell(w,\nu,\lambda)&\lesssim \int_{\S^2}\frac{1}{(1-\lambda|\hat{w}\cdot\sigma|)^\ell}\,d\sigma= 4\pi\int_0^1\frac{1}{(1-\lambda x)^{\ell/2}}\,dx\\
&=\frac{4\pi}{\lambda}\int_0^\lambda \frac{1}{(1-z)^{\ell}}\,dz\lesssim  \frac{1}{\lambda(\ell-1)}(1-\lambda)^{1-\ell},
\end{align*}
since $\ell>1$.   
\end{proof}

\begin{lemma}\label{lem: perp lemma}
Let $\ell>1$.  For $w,\nu\in\R^3\setminus\{0\}$, we define 
 \begin{equation}\label{definition of B}
\mathcal{B}_{\ell}(w,\nu):= \mathds{1}_{|\hat{w}\cdot\hat{\nu}|<\sqrt{2}/2}\int_{\S^2}\left(\la w\ra^2|\hat{\nu}\cdot\sigma|^2-2|w||\nu||\hat{w}\cdot\sigma||\hat{\nu}\cdot\sigma|+|\nu|^2\right)^{-\ell}\,d\sigma.     
 \end{equation}
 Then, there holds the estimate
 \begin{equation}\label{B estimate}
 \mathcal{B}_{\ell}(w,\nu)     \lesssim \frac{\la w\ra^{2\ell-2}}{|w|}|\nu|^{1-2\ell}
 \end{equation}  
\end{lemma}
\begin{proof} 
Let $w,\nu\neq 0$ with $|\hat{w}\cdot\hat{\nu}|<\sqrt{2}/2$. We parametrize the sphere by choosing the  $z$-axis in the direction of $\nu$ and the $x-$axis along the projection of $w$ on $\nu^\perp$. Let $\theta\in[0,\pi]$ be  the polar angle and $\phi\in[0,2\pi)$ the azimuthal angle. Then, any $\sigma\in\S^2$ can be uniquely written as $\sigma=(\cos\phi\sin\theta,\sin\phi\sin\theta,\cos\theta)$. In particular, we have  $\hat{\nu}\cdot\sigma=\cos\theta$ and $\hat{w}=(\sin\psi,0,\cos\psi)$, where $\psi=\arccos(\hat{w}\cdot\hat{\nu})$, therefore
$$\hat{w}\cdot\sigma=\cos\phi\sin\theta\sin\psi+\cos\theta\cos\psi.$$
Using the Cauchy-Schwarz inequality in $\R^2$, we estimate
\begin{align*}
|\hat{w}\cdot\sigma|&=|\cos\phi\sin\theta\sin\psi+\cos\theta\cos\psi|=|(\cos\phi\sin\psi,\cos\psi)\cdot(\sin\theta,\cos\theta)|\notag\\
&\leq \sqrt{\cos^2\phi\sin^2\psi+\cos^2\psi}=\sqrt{1-\sin^2\phi\sin^2\psi},
\end{align*}
therefore
\begin{equation}\label{bound by mu}
|\hat{w}\cdot\sigma|\mathds{1}_{|\hat{w}\cdot\hat{\nu}|<\sqrt{2}/2}\leq \mathds{1}_{\sin\psi>\sqrt{2}/2}\sqrt{1-\sin^2\phi\sin^2\psi}<\sqrt{1-\frac{1}{2}\sin^2\phi}:=\mu_\phi.    
\end{equation}
The crucial fact is that the above bound is uniform in $\theta$.

Namely, integrating in spherical coordinates and using \eqref{bound by mu} we obtain
\begin{align*}
\mathcal{B}(w,\nu)
&\leq \int_0^{2\pi}\int_{0}^\pi\frac{\sin\theta}{\left(\l w\r^2\cos^2\theta-2\mu_\phi|w||\nu||\cos\theta|+|\nu|^2\right)^{\ell}}\,d\theta\,d\phi\\
&\approx \int_0^{2\pi}\int_0^1\frac{1}{\left(\l w\r^2x^2-2\mu_\phi|w||\nu|x+|\nu|^2\right)^{\ell}}\,dx\,d\phi\\
&=\int_0^{2\pi}\int_0^1 P^{-\ell}_{w,\nu,\phi}(x)\,dx\,d\phi,
\end{align*}
where 
$$P_{w,\nu,\phi}(x):=\la w\ra^2 x^2-2\mu_\phi |w||\nu| x+|\nu|^2.$$
We  factor
$$P_{w,\nu,\phi}(x)=\l w\r^{2}b_{w,\nu,\phi}^{2}\left(\left(\frac{x-x_{w,\nu,\phi}}{b_{w,\nu,\phi}}\right)^2+1\right),$$
where 
\begin{align*}
x_{w,\nu,\phi}&:=\frac{|w||\nu|\mu_\phi}{\l w\r^2},\qquad b_{w,\nu,\phi}:=\frac{|\nu|}{\l w\r^2}\left(1+|w|^2\left(1-\mu^2_\phi\right)\right)^{1/2}.
\end{align*}
Using the definition of $\mu_\phi$ from \eqref{bound by mu}, we obtain
\begin{equation}\label{b without mu}
b_{w,\nu,\phi}=  \frac{|\nu|}{\l w\r^2}\left(1+\frac{|w|^2}{2}\sin^2\phi\right)^{1/2}  
\end{equation}

Therefore,
\begin{align}
\mathcal{B}(w,\nu)&\lesssim \l w\r^{-2\ell}\int_0^{2\pi}b_{w,\nu,\phi}^{-2\ell}\int_0^1 \left(\left(\frac{x-x_{w,\nu,\phi}}{b_{w,\nu,\phi}}\right)^2+1\right)^{-\ell}\,dx\,d\phi\notag\\
&\leq \l w\r^{-2\ell}\int_0^{2\pi}b_{w,\nu,\phi}^{1-2\ell}\int_{-\infty}^{+\infty}(1+z^2)^{-\ell}\,dz\notag\\
&\lesssim \la w\ra^{-2\ell}\int_0^{2\pi}b_{w,\nu,\phi}^{1-2\ell}\,d\phi\notag\\
&= \l w\r^{2\ell-2}|\nu|^{1-2\ell}\int_0^{2\pi}\left(1+\frac{|w|^2}{2}\sin^2\phi\right)^{1/2-\ell}\,d\phi,\label{first estimate on B}
\end{align}
where for the second to last inequality we used the fact $\ell>1$ for the convergence of the integral in $z$, and for the last equality we used \eqref{b without mu}. 

Now the function $\phi\mapsto \left(1+\frac{|w|^2}{2}\sin^2\phi\right)^{1/2-\ell}$ is  $\pi$-periodic and even, thus we have
\begin{align}
 \int_0^{2\pi}(1+|w|^2\sin^2\phi)^{1/2-\ell}\,d\phi&=4\int_0^{\pi/2} (1+|w|^2\sin^2\phi)^{1/2-\ell}\,d\phi\approx \int_0^{\pi/2}(1+|w|^2\phi^2)^{1/2-\ell}\,d\phi\notag\\
 &\lesssim \frac{1}{|w|}\int_0^\infty (1+z^2)^{1/2-\ell}\,dz\lesssim \frac{1}{|w|},\label{second bound on B}
\end{align}
where we used the facts $\sin\phi\approx\phi$, $\phi\in[0,\pi/2]$ and  $\ell>1$ for the convergence of the integral in $z$. 

Combining \eqref{first estimate on B}, \eqref{second bound on B}, estimate \eqref{B estimate} follows.

\end{proof}

\subsection{Decomposition of the hard-potential operator} Let $f,g\in C_c(\R^3)$. Since $b\in L^\infty([-1,1])$, we 
can bound
\begin{equation}\label{L infty bilinear estimate}
   |Q^+(f,g)(v)|\leq \|b\|_{L^\infty}\int_{\R^3\times\S^2}|u|^\gamma|f(v^*)g(v_1^*)|\,d\sigma\,dv_1. 
\end{equation}
Then, we
decompose the operator into small, regular and singular pieces as follows:
\begin{align}
|Q^+(f,g)(v)|
&\lesssim \mathcal{K}_{sm}(f,g)(v)+\mathcal{K}_{reg}(f,g)(v)+\mathcal{K}_{sin}(f,g)(v), \label{decomposition of Q},   
\end{align}
where
\begin{align}
\K_{sm}(f,g)(v)&=\int_{\R^3\times\S^2}|u|^\gamma |f(v^*)g(v_1^*)|\mathds{1}_{|v|\leq 3}\mathds{1}_{|v_1^*|^2> E/2}\,d\sigma\,dv_1,\label{K small}\\
\K_{reg}(f,g)(v)&=\int_{\R^3\times\S^2}|u|^\gamma |f(v^*)g(v_1^*)|\mathds{1}_{|v|>3}\mathds{1}_{|v_1^*|^2> E/2}\mathds{1}_{|v^*|\geq |v|/6}\,d\sigma\,dv_1,\label{K reg}\\
\K_{sin}(f,g)(v)&=\int_{\R^3\times\S^2}|u|^\gamma  |f(v^*)g(v_1^*)|\mathds{1}_{|v|>3}\mathds{1}_{|v_1^*|^2> E/2}\mathds{1}_{|v^*|<|v|/6}\,d\sigma\,dv_1.\label{K sin}
\end{align}

\subsection{The $L_v^\infty$-- estimate}
\begin{proposition}\label{Q1 proposition Linf}
Let $0<\gamma <2 $ and $b$ satisfying \eqref{assumption 1}--\eqref{assumption 3}. Let $1< p,q< 3$ with $p\geq \frac{2}{2-\gamma}$, such that 
$\frac{1}{p}+\frac{1}{q}=1$, and fix $k\geq l> \frac{2}{q}$. Then, there holds the estimate
\begin{equation}\label{L inf bound Q_1}
\|\l v\r^{k} Q^+(f,g)\|_{L_v^\infty}\lesssim \|\l v\r^l f\|_{L_v^{p}}\|\l v\r^k g\|_{L_v^{q}},\quad\forall\,f,g\in C_c(\R^3).
\end{equation}
\end{proposition}
\begin{proof} 
Let  $f,g\in C_c(\R^3)$, and denote $f_l:=\l v\r^l f$, $g_k:=\l v\r^k g$. By \eqref{decomposition of Q}, it suffices to bound $\mathcal{K}_{sm}$, $\mathcal{K}_{reg}$ and $\mathcal{K}_{sin}$.

Before we start, we note $p\geq \frac{2}{2-\gamma}$ is equivalent to $q\leq 2/\gamma$. Since $ k\geq l> 2/q $, we have $k\geq l>\gamma$. Moreover, since $q<3$, we have $q<2p$ thus $k\geq l> 2/q>1/p$.

\medskip
\noindent {\bf Estimate for $\K_{sm},\, \mathcal{K}_{reg}$}  Using the inequality $|u|\leq |v|+|v_1|\lesssim E^{1/2}$, we bound 
\begin{align}
\frac{|u|^\gamma\l v\r^{k}}{\l v^*\r^l \l v_1^*\r^k}\mathds{1}_{|v|\leq 4}\mathds{1}_{|v_1^*|^2> E/2}&\lesssim (1+E)^{\frac{\gamma-k}{2}}\lesssim 1,\label{pointwise bound small}\\
\frac{|u|^\gamma \l v\r^{k}}{\l v^*\r^l\l v_1^*\r^k}\mathds{1}_{|v|>3}\mathds{1}_{|v^*|\geq|v|/6}\mathds{1}_{|v_1^*|^2>E/2}&\lesssim (1+E)^{\frac{\gamma}{2}-\frac{k}{2}} \la v\r^{k-l}\leq \la v\r^{\gamma-l}\lesssim 1\label{pointwise bound reg}
\end{align}
since  $k\geq l> \gamma$.

Recalling the Maxwell operator notation \eqref{Maxwell notation}, and using the above estimates, we obtain
\begin{align}
\l v\r^{k} \K_{sm}(f,g)(v)
&\lesssim \int_{\R^3\times\S^2} |f_l(v^*) g_k(v_1^*)|\,d\sigma\,dv_1=Q_{M,1}^+(|f_l|,|g_k|), \label{K bound small}\\
\l v\r^{k} \K_{sm}(f,g)(v)
&\lesssim \int_{\R^3\times\S^2} |f_l(v^*) g_k(v_1^*)|\,d\sigma\,dv_1=Q_{M,1}^+(|f_l|,|g_k|).\label{K bound regular}
\end{align}
Since $1<p,q<3$, Theorem \ref{Maxwell theorem general} implies
\begin{align}
\|\l v\r^{k} \K_{sm}(f,g)\|_{L_v^\infty}&\lesssim \|f_l\|_{L_v^{p}}\|g_k\|_{L_v^{q}},  \label{Linf bound Q1 sm}\\
\|\l v\r^{k} \K_{reg}(f,g)\|_{L_v^\infty}&\lesssim \|f_l\|_{L_v^{p}}\|g_k\|_{L_v^{q}}.\label{Linf bound Q1 reg}
\end{align}

\medskip
\noindent It remains to estimate the singular operator $\K_{sin}$.  

\medskip
\noindent{\bf Estimate for $\K_{sin}$} We start with a geometric implications of the fact that $|v^*|<|v|/6$ in the corresponding support. Namely,using the triangle inequality and \eqref{R+ v*}, \eqref{magnitude}, we have 
\begin{equation}\label{upper bound on v}
|v|\leq |v-v^*|+|v^*|\leq |R_\sigma^+(u)|+\frac{|v|}{6}\Rightarrow |v|\leq \frac{6}{5}|R_\sigma^+(u)|=\frac{6}{5}|u||\widehat{R}_\sigma^+(u)\cdot\sigma|.
\end{equation}

\noindent Now the above bound followed by the inequality $|u|\lesssim E^{1/2}$ imply
\begin{align}
  \frac{ \l v\r^{k}}{ \l v_1^*\r^k}&\mathds{1}_{|v|>3}\mathds{1}_{|v_1^*|^2>E/2}\mathds{1}_{|v^*|<|v|/6}\lesssim \frac{ |v|^k}{\la v_1^*\ra^k}\mathds{1}_{|v|>3}\mathds{1}_{|v_1^*|^2>E/2}\mathds{1}_{|v^*|<|v|/6}\notag\\
  &\lesssim \frac{|u|^{k}|\hat{R}^+_\sigma(u)\cdot\sigma|^k}{\la v_1^*\ra^k}\mathds{1}_{|v|>3}\mathds{1}_{|v_1^*|^2>E/2}\mathds{1}_{|v^*|<|v|/6}\notag\\
  &\lesssim |\hat{R}^{+}_\sigma(u)\cdot\sigma|^k\mathds{1}_{|v|>3}\mathds{1}_{|v_1^*|^2>E/2}\mathds{1}_{|v^*|<|v|/6}\label{pointwise cancellation bound}.
\end{align}

 Since $p<3$, we have $p<2q$. Pick $1<\alpha<2q/p$. Using \eqref{pointwise cancellation bound} and H\"older's inequality, we have
\begin{align}
\la v\ra^k& \K_{sin}(f,g)(v)=\int_{\R^3\times\S^2} \frac{|u|^\gamma \l v\r^{k}}{\l v^*\r^l \l v_1^*\r^k}\mathds{1}_{|v|>3}\mathds{1}_{|v_1^*|^2>E/2}\mathds{1}_{|v^*|<|v|/6} |f_l(v^*)g_k(v_1^*)|\,d\sigma\,dv_1\notag\\
&\lesssim \int_{\R^3\times\S^2}\left(|\hat{R}_\sigma^+(u)\cdot\sigma|^k |f_l(v^*)|\right)\left(\frac{|u|^\gamma}{\la v^*\ra^l}|g_k(v_1^*)|\mathds{1}_{|v|>3}\mathds{1}_{|v_1^*|^2>E/2}\mathds{1}_{|\hat{v_1^*}\cdot\sigma|>\sqrt{3}/2}\right)\,d\sigma\,dv_1\notag\\
&=\int_{\R^3\times\S^2}\left(\frac{|\hat{R}_\sigma^+(u)\cdot\sigma|^k}{|\hat{R}_\sigma^-(u)\cdot\sigma|^{\alpha/q}}|f_l(v^*)|\right)\left(\frac{|u|^\gamma |\hat{R}_\sigma^-(u)\cdot\sigma|^{\alpha/q}}{\la v^*\ra^l}|g_k(v_1^*)|\mathds{1}_{|v|>3}\mathds{1}_{|v_1^*|^2>E/2}\mathds{1}_{|\hat{v_1^*}\cdot\sigma|>\sqrt{3}/2}\right)\,d\sigma\,dv_1\notag\\
&\leq \mathcal{J}_1^{1/p}\mathcal{J}_2^{1/q},\label{Holder nonperp Linfty}
\end{align}
where 
\begin{align}
\mathcal{J}_1&=\int_{\R^3\times\S^2}\frac{|\hat{R}_\sigma^+(u)\cdot\sigma|^{kp}}{|\hat{R}_\sigma^-(u)\cdot\sigma|^{\alpha p/q}}|f_l(v^*)|^p\,d\sigma\,dv_1\label{I1 nonperp Linfty}\\
\mathcal{J}_2&=\int_{\R^3\times\S^2}\frac{|u|^{\gamma q}|\hat{R}_\sigma^-(u)\cdot\sigma|^\alpha}{\la v^*\ra^{l q}}|g_k(v_1^*)|^q\mathds{1}_{|v|>3}\mathds{1}_{|v_1^*|^2>E/2}\mathds{1}_{|v^*|<|v|/6}\,d\sigma\,dv_1.\label{I2 nonperp Linfty}
\end{align}
Since $kp>1$ and $\alpha p/q<2$, estimate \eqref{averaging estimate on I_2 Maxwell} implies
\begin{equation}\label{J1 bound Linfty}
\mathcal{J}_1\lesssim \|f_l\|_{L^p}^p.
\end{equation}

\noindent For $\mathcal{J}_2$, we decompose in two regions based on whether $v_1^*$ and $R_\sigma^-(u)=v-v_1^*$ are closer to aligned or perpendicular. Namely, we write 
\begin{equation}\label{decomposition of I2}
\mathcal{J}_2=\mathcal{J}_2^{alig}+\mathcal{J}_2^{perp},    
\end{equation}
where
\begin{align}
\mathcal J_2^{alig}&= \int_{\R^3\times\S^2}\frac{|u|^{\gamma q}|\hat{R}_\sigma^-(u)\cdot\sigma|^\alpha}{\la v^*\ra^{l q}}|g_k(v_1^*)|^q\mathds{1}_{|v|>3}\mathds{1}_{|v_1^*|^2>E/2}\mathds{1}_{|v^*|<|v|/6}\mathds{1}_{|\hat{v_1^*}\cdot\hat{R}_\sigma^-(u)|\geq \sqrt{2}/2}\,d\sigma\,dv_1 \\
\mathcal J_2^{perp}&= \int_{\R^3\times\S^2}\frac{|u|^{\gamma q}|\hat{R}_\sigma^-(u)\cdot\sigma|^\alpha}{\la v^*\ra^{l q}}|g_k(v_1^*)|^q\mathds{1}_{|v|>3}\mathds{1}_{|v_1^*|^2>E/2}\mathds{1}_{|v^*|<|v|/6}\mathds{1}_{|\hat{v_1^*}\cdot\hat{R}_\sigma^-(u)|< \sqrt{2}/2}\,d\sigma\,dv_1 
\end{align}

\subsubsection*{Estimate for $\mathcal{J}_2^{alig}$} First, we note that in the corresponding support we have $|\hat{v_1^*}\cdot\sigma|>\sqrt{3}/2$. Indeed, if $|\hat{v_1^*}\cdot\sigma|\leq \sqrt{3}/2$, estimate \eqref{lower bound on post-velocity magnitude} and the fact that $|v_1^*|^2>E/2$ imply
$$|v^*|^2\geq (1-\frac{\sqrt{3}}{2})|v_1^*|^2>\frac{|v_1^*|^2}{14}>\frac{E}{28}>\frac{|v|^2}{36}\Rightarrow |v^*|>\frac{|v|}{6}$$

\noindent For $a\in\R^{3}\setminus \{0\}$,  let us denote $\zeta_a:=\frac{|a|^2}{\la a\ra^2}$.
Using \eqref{lower bound for nonperp} and the inequality $|u|\lesssim E^{1/2}$, we have
\begin{align}
\frac{|u|^\gamma}{\la v^*\ra^l}\mathds{1}_{|v|>3}\mathds{1}_{|v_1^*|^2>E/2}&\lesssim \frac{|u|^\gamma}{\la v_1^*\ra^l\left(1-\zeta_{v_1^*}|\hat{v_1^*}\cdot\sigma|\right)^{l/2}}\mathds{1}_{|v|>3}\mathds{1}_{|v_1^*|^2>E/2}\notag\\
&\lesssim  \frac{\mathds{1}_{|v_1^*|>2}}{\la v_1^*\ra^{l-\gamma}\left(1-\zeta_{v_1^*}|\hat{v_1^*}\cdot\sigma|\right)^{l/2}},\label{cancellation bound alig}    
\end{align}
where for the last inequality we used the fact that $\mathds{1}_{|v|>3}\mathds{1}_{|v_1^*|^2>E/2}\leq \mathds{1}_{|v_1^*|>2}$ since $|v_1^*|^2>E/2>|v|^2/2>4$.

Using \eqref{cancellation bound alig}, followed by \eqref{R- v_1*} and the substitution $y:=u=v-v_1$, we obtain
\begin{align*}
\mathcal{J}_2^{alig}&\lesssim  \int_{\R^3\times\S^2}\frac{|g_k(v-R_\sigma^-(u))|^q\mathds{1}_{|v-R_\sigma^-(u)|>2}\mathds{1}_{|\hat{(v-R_\sigma^-(u))}\cdot\sigma|>\sqrt{3}/2}\mathds{1}_{|\hat{(v-R_\sigma^-(u))}\cdot \hat{R}_\sigma^-(u)|\geq\sqrt{2}/2}}{\la v-R_\sigma^-(u)\ra^{(l-\gamma)q}\left(1-\zeta_{v-R_\sigma^-(u)}|\hat{(v-R_\sigma(u))}\cdot\sigma|\right)^{lq/2}}\,d\sigma\,dv_1 \\
&=\int_{\S^2}\int_{\R^3}\frac{|g_k(v-R_\sigma^-(y))|^q\mathds{1}_{|v-R_\sigma^-(y)|>2}\mathds{1}_{|\hat{(v-R_\sigma^-(y))}\cdot\sigma|>\sqrt{3}/2}\mathds{1}_{|\hat{(v-R_\sigma^-(y))}\cdot \hat{R}_\sigma^-(y)|\geq\sqrt{2}/2}}{\la v-R_\sigma^-(y)\ra^{(l-\gamma)q}\left(1-\zeta_{v-R_\sigma^-(y)}|\hat{(v-R_\sigma(y))}\cdot\sigma|\right)^{lq/2}}\,dy\,d\sigma.
\end{align*}
For fixed $\sigma\in\S^2$, we use Proposition \ref{chgvarkin} to change variables $\nu:=R_\sigma^-(y)$, and obtain
\begin{align*}
I_2&\approx \int_{\S^2}\int_{\R^3}\frac{|g_k(v-\nu)|^q\mathds{1}_{|v-\nu|>2}\mathds{1}_{|\hat{(v-\nu)}\cdot\sigma|>\sqrt{3}/2}\mathds{1}_{|\hat{(v-\nu)}\cdot\sigma|>\sqrt{2}/2}}{\la v-\nu\ra^{(l-\gamma)q}|\hat{\nu}\cdot\sigma|^2\left(1-\zeta_{v-\nu}|\hat{v-\nu}\cdot\sigma|\right)^{lq/2}}\,d\nu\,d\sigma \\
&=\int_{\R^3}\frac{|g_k(v-\nu)|^q\mathds{1}_{|v-\nu|>2}}{\la v-\nu\ra^{(l-\gamma)q}}\left(\int_{\S^2}\frac{\mathds{1}_{|\hat{(v-\nu)}\cdot\sigma|>\sqrt{3}/2}\mathds{1}_{|\hat{(v-\nu)}\cdot\sigma|>\sqrt{2}/2}}{|\hat{\nu}\cdot\sigma|^2\left(1-\zeta_{v-\nu}|\hat{v-\nu}\cdot\sigma|\right)^{lq/2}}\,d\sigma\right)\,d\nu\\
&=\int_{\R^3}\frac{|g_k(v-\nu)|^q\mathds{1}_{|v-\nu|>2}}{\la v-\nu\ra^{(l-\gamma)q}}\mathcal{A}_{lq/2}(v-\nu,\nu,\zeta_{v-\nu})\,d\nu,
\end{align*}
where we use notation from \eqref{definition of A}. Since $l>2/q$, Lemma \ref{lem: lemma on A} implies 
\begin{align*}\mathcal{A}_{lq/2}(v-\nu,\nu,\zeta_{v-\nu})\lesssim \frac{1}{\zeta_{v-\nu}}(1-\zeta_{v-\nu})^{1-ql/2}=\frac{\la v-\nu\ra^2}{|v-\nu|^2}\la v-\nu\ra^{ql-2},
\end{align*}
thus,  since $q\leq 2/\gamma$, we obtain 
\begin{align}\label{estimate on J alig}
\mathcal{J}_2^{alig}\lesssim\int_{\R^3}\mathds{1}_{|v-\nu|>2}\frac{\la v-\nu\ra^2}{|v-\nu|^2}\la v-\nu\ra^{q\gamma-2} |g_k(v-\nu)|^q\,d\nu\lesssim \int_{\R^3}|g_k(v-\nu)|^q\,d\nu=\|g_k\|_{L^q}^q,  
\end{align}
by translation invariance.

\subsubsection*{Estimate for $\mathcal{J}^{perp}$}
Using \eqref{lower bound for perp} and \eqref{magnitude}, we have 
\begin{align*}&\frac{|u|^{\gamma q} |\hat{R}_\sigma^-(u)\cdot\sigma|^\alpha}{\la v^*\ra^l}\leq \frac{|u|^{\gamma q}|\hat{R}_\sigma^-(u)|^{l+\alpha}}{\left(\la v_1^*\ra^2|\hat{R}_\sigma^-(u)\cdot\sigma|^2-2|v_1^*||R_\sigma^-(u)||\hat{v_1^*}\cdot\sigma||\hat{R}_\sigma^-(u)\cdot\sigma|+|R_\sigma^-(u)|^2\right)^{l/2}}\\
& =\frac{|u|^{\gamma q+2-l-\alpha}|R_\sigma^-(u)|^{l+\alpha-2}|\hat{R}_\sigma^-(u)\cdot\sigma|^2}{\left(\la v_1^*\ra^2|\hat{R}_\sigma^-(u)\cdot\sigma|^2-2|v_1^*||R_\sigma^-(u)||\hat{v_1^*}\cdot\sigma||\hat{R}_\sigma^-(u)\cdot\sigma|+|R_\sigma^-(u)|^2\right)^{l/2}}
\end{align*}
Now for $|v|>3$, $|v_1^*|^2>E/2$ and $|v^*|<|v|/6$, we have $|u|\approx |v_1^*|$, $|R_\sigma^-(u)|\leq 3|v_1^*|$ and $|v_1^*|>2$. Indeed, 
\begin{align*}
|u|&\lesssim E^{1/2}\lesssim |v_1^*|,\quad |u|=|v^*-v_1^*|\geq |v_1^*|-|v^*|>|v_1^*|-\frac{|v|}{6}>(1-\frac{\sqrt{2}}{6})|v_1^*|,\\
|R_\sigma^-(u)|&=|v-v_1^*|\leq |v|+|v_1^*|\leq E^{1/2}+|v_1^*|\leq 3|v_1^*|,\quad 
|v_1^*|^2>E/2>|v^*|^2>2>4.
\end{align*}
Therefore, we can bound 
\begin{align*}
&\frac{|u|^{\gamma q} |\hat{R}_\sigma^-(u)\cdot\sigma|^\alpha}{\la v^*\ra^{lq}}\mathds{1}_{|v|>3}\mathds{1}_{|v_1^*|^2>E/2}\mathds{1}_{|v^*|<|v|/6}\\
&\lesssim \frac{|v_1^*|^{\gamma q+2-lq-\alpha}|R_\sigma^-(u)|^{lq+\alpha-2}|\hat{R}_\sigma^-(u)\cdot\sigma|^2}{\left(\la v_1^*\ra^2|\hat{R}_\sigma^-(u)\cdot\sigma|^2-2|v_1^*||R_\sigma^-(u)||\hat{v_1^*}\cdot\sigma||\hat{R}_\sigma^-(u)\cdot\sigma|+|R_\sigma^-(u)|^2\right)^{lq/2}}\mathds{1}_{|v_1^*|>3}\mathds{1}_{|R_\sigma^-(u)|\leq 3|v_1^*|}.
\end{align*}
Using the above bound followed by \eqref{R- v_1*} and the substitution $y:=u=v-v_1$, we obtain 
\begin{align*}
&\mathcal{J}_2^{perp}\lesssim\int_{\R^3\times\S^2}\frac{|v_1^*|^{\gamma q+2-lq-\alpha}|R_\sigma^-(u)|^{lq+\alpha-2}|\hat{R}_\sigma^-(u)\cdot\sigma|^2|g_k(v_1^*)|^q }{\left(\la v_1^*\ra^2|\hat{R}_\sigma^-(u)\cdot\sigma|^2-2|v_1^*||R_\sigma^-(u)||\hat{v_1^*}\cdot\sigma||\hat{R}_\sigma^-(u)\cdot\sigma|+|R_\sigma^-(u)|^2\right)^{lq/2}}\\
&\hspace{5cm}\times\mathds{1}_{|v_1^*|>2}\mathds{1}_{|R_\sigma^-(u)|\leq 3|v_1^*|}\mathds{1}_{|\hat{v_1^*}\cdot\hat{R}_\sigma^-(u)|<\sqrt{2}/2}\,d\sigma\,dv_1\\
&=\int_{\R^3\times\S^2}\frac{|v-R_\sigma^-(u)|^{\gamma q+2-lq-\alpha}|R_\sigma^-(u)|^{lq+\alpha-2}|\hat{R}_\sigma^-(u)\cdot\sigma|^2|g_k(v-R_\sigma^-(u))|^q }{\left(\la v-R_\sigma^-(u)\ra^2|\hat{R}_\sigma^-(u)\cdot\sigma|^2-2|v-R_\sigma^-(u)||R_\sigma^-(u)||\hat{(v-R_\sigma^-(u))}\cdot\sigma||\hat{R}_\sigma^-(u)\cdot\sigma|+|R_\sigma^-(u)|^2\right)^{lq/2}}\\
&\hspace{5cm}\times \mathds{1}_{|v-R_\sigma^-(u)|>2}\mathds{1}_{|R_\sigma^-(u)\leq 3|v-R_\sigma^-(u)|}\mathds{1}_{|\hat{(v-R_\sigma^-(u))}\cdot\hat{R}_\sigma^-(u)|<\sqrt{2}/2}\,d\sigma\,dv_1\\
&=\int_{\S^2}\int_{\R^3}\frac{|v-R_\sigma^-(y)|^{\gamma q+2-lq-\alpha}|R_\sigma^-(y)|^{lq+\alpha-2}|\hat{R}_\sigma^-(y)\cdot\sigma|^2|g_k(v-R_\sigma^-(y))|^q\mathds{1}_{|\hat{(v-R_\sigma^-(y))}\cdot\hat{R}_\sigma^-(y)|<\sqrt{2}/2}}{\left(\la v-R_\sigma^-(y)\ra^2|\hat{R}_\sigma^-(y)\cdot\sigma|^2-2|v-R_\sigma^-(y)||R_\sigma^-(y)||\hat{(v-R_\sigma^-(y))}\cdot\sigma||\hat{R}_\sigma^-(y)\cdot\sigma|+|R_\sigma^-(y)|^2\right)^{lq/2}}\\
&\hspace{5cm}\times \mathds{1}_{|v-R_\sigma^-(y)|>2}\mathds{1}_{|R_\sigma^-(y)\leq 3|v-R_\sigma^-(y)|}\mathds{1}_{|\hat{(v-R_\sigma^-(y))}\cdot\hat{R}_\sigma^-(y)|<\sqrt{2}/2}\,dy\,d\sigma\\
\end{align*}
For fixed $\sigma\in\S^2$, we use Proposition \ref{chgvarkin} to substitute $\nu:=R_\sigma^-(y)$ and obtain
\begin{align*}
\mathcal{J}_2^{perp}&\lesssim \int_{\S^2}\int_{\R^3}\frac{|v-\nu|^{\gamma q+2-lq-\alpha}|\nu|^{lq+\alpha-2}|g_k(v-\nu)|^q\mathds{1}_{|v-\nu|>3}\mathds{1}_{|\nu|\leq 3|v-\nu|}\mathds{1}_{|\hat{(v-\nu)}\cdot\hat{\nu}|<\sqrt{2}/2} }{\left(\la v-\nu\ra^2|\hat{\nu}\cdot\sigma|^2-2|v-\nu||\nu||\hat{(v-\nu)}\cdot\sigma||\hat{\nu}\cdot\sigma|+|\nu|^2\right)^{lq/2}}  \,d\nu\,d\sigma\\
&=\int_{\R^3}|v-\nu|^{\gamma q+2-lq-\alpha}|\nu|^{lq+\alpha-2}\mathds{1}_{|v-\nu|>3}\mathds{1}_{|\nu|\leq 3|v-\nu|}\mathds{1}_{|\hat{(v-\nu)}\cdot\hat{\nu}|<\sqrt{2}/2}\\
&\quad\times\int_{\S^2}\left(| v-\nu|^2|\hat{\nu}\cdot\sigma|^2-2|v-\nu||\nu||\hat{(v-\nu)}\cdot\sigma||\hat{\nu}\cdot\sigma|+|\nu|^2\right)^{-lq/2}\,d\sigma\,d\nu\\
&=\int_{\R^3}| v-\nu|^{\gamma q+2-lq-\alpha}|\nu|^{lq+\alpha-2}\mathcal{B}_{lq/2}(v-\nu,\nu)\mathds{1}_{|v-\nu|>3}\mathds{1}_{|\nu|\leq 3|v-\nu|}\,d\nu,
\end{align*}
where we recall ${\mathcal{B}_{lq/2}}$ from \eqref{definition of B}. Since $l>2/q$, Lemma \ref{lem: lemma on A} implies 
$$\mathcal{B}_{lq/2}\lesssim \frac{\la v-\nu\ra^{lq-2}}{|v-\nu|}|\nu|^{1-lq},$$
thus since $\alpha>1$ and $q\leq 2/\gamma$, we obtain
\begin{align}
\mathcal{J}_2^{perp}&\lesssim \int_{\R^3}| v-\nu|^{\gamma q+1+lq-\alpha}\la v-\nu\ra^{lq-2}|\nu|^{\alpha-1}|g_k(v-\nu)|^q\mathds{1}_{|v-\nu|>3}\mathds{1}_{|\nu|\leq 3|v-\nu|}\,d\nu\notag\\
&\lesssim \int_{\R^3}\la v-\nu\ra^{\gamma q -2}|g_k(v-\nu)|^q\,d\nu\leq \int_{\R^3}|g_k(v-\nu)|^q\,d\nu=\|g_k\|_{L^q}^q.\label{estimate on Jperp}
\end{align}

Combining \eqref{estimate on J alig}, \eqref{estimate on Jperp}, \eqref{decomposition of I2}, we obtain 
\begin{equation}\label{J2 bound Linfty}
\mathcal{J}_2\lesssim \|g_k\|_{L^q}^q.    
\end{equation}
Finally, putting together \eqref{Holder nonperp Linfty}, \eqref{J1 bound Linfty}, \eqref{J2 bound Linfty}, we conclude 
\begin{equation}
\|\la v\ra^k\mathcal{K}_{sin}(f,g)\|_{L^\infty}\lesssim \|f_l\|_{L^p}\|g\|_{L^q},     
\end{equation}
which combined with \eqref{Linf bound Q1 sm}, \eqref{Linf bound Q1 reg} and \eqref{decomposition of Q} yields \eqref{L inf bound Q_1}.
The proof is complete.
\end{proof}

\subsection{The $L_v^r$-- estimate}

\begin{proposition}\label{Q_1 proposition Lr}
Let $0<\gamma <2 $ and $b$ satisfying \eqref{assumption 1}--\eqref{assumption 3}. Let $1\leq p,q,<3$ with $p\geq \frac{2}{2-\gamma}$, and $\frac{2}{2-\gamma}\leq r< \infty$, such that 
$\frac{1}{p}+\frac{1}{q}=1+\frac{1}{r}$, and fix $k\geq l> 2/p'$ with $k> 1/q'$. Then, there holds the estimate
\begin{equation}\label{L inf bound Q_1}
\|\l v\r^{k} Q^+(f,g)\|_{L_v^r}\lesssim \|\l v\r^l f\|_{L_v^{p}}\|\l v\r^k g\|_{L_v^{q}},\qquad \forall\, f,g\in C_c(\R^3).
\end{equation}
\end{proposition}
\begin{proof}
We use again decomposition \eqref{decomposition of Q}. Again, it suffices to estimate $\mathcal{K}_{sm},\,\mathcal{K}_{reg},\,\mathcal{K}_{sin}$.

We note that the assumption $p\geq \frac{2}{2-\gamma}$ is equivalent to $p'\leq 2/\gamma$. Thus we obtain $k\geq l>2/p'\geq \gamma$. 

\medskip
\noindent{\bf Estimate for $\mathcal{K}_{sm},\mathcal{K}_{reg}$} Since $k\geq l>\gamma$ \eqref{pointwise bound small}--\eqref{K bound regular} hold,  thus we obtain
\begin{align*}
\|\la v\r^k \mathcal{K}_{sm}(f,g)\|_{L^r}&\lesssim \|Q_{M,1}^+(|f_l|,|g_k|)\|_{L^r}\\ 
\|\la v\r^k \mathcal{K}_{reg}(f,g)\|_{L^r}&\lesssim \|Q_{M,1}^+(|f_l|,|g_k|)\|_{L^r}.
\end{align*}
Since $1\leq p,q<3$, Theorem \ref{Maxwell theorem general} implies
\begin{align}
  \|\la v\r^k \mathcal{K}_{sm}(f,g)\|_{L^r}&\lesssim \|f_l\|_{L^p}\|g\|_{L^q}\label{estimate Ksm Lr},\\ 
\|\la v\r^k \mathcal{K}_{reg}(f,g)\|_{L^r}&\lesssim \|f_l\|_{L^p}\|g\|_{L^q}\label{estimate Kreg Lr}.  
\end{align}

\medskip
\noindent{\bf Estimate for $\mathcal{K}_{sin}$} Since $p<3$, the same argument as in the proof of Proposition \ref{Maxwell prop Lr} yields existence of $\alpha>1$ and $s\in \R$ with 
\begin{equation}\label{s condition hard}
\frac{(1-s)\alpha q'}{p'}<2,\qquad\frac{s\alpha r}{p'}<2.
\end{equation}

Using \eqref{pointwise cancellation bound}, the exponents' relation and H\"older's inequality, we obtain
\begin{align}
\la v\ra^k& \K_{sin}(f,g)(v)=\int_{\R^3\times\S^2} \frac{|u|^\gamma \l v\r^{k}}{\l v^*\r^l \l v_1^*\r^k}\mathds{1}_{|v|>3}\mathds{1}_{|v_1^*|^2>E/2}\mathds{1}_{|v^*|<|v|/6} |f_l(v^*)g_k(v_1^*)|\,d\sigma\,dv_1\notag\\
&\lesssim \int_{\R^3\times\S^2}\left(|\hat{R}_\sigma^+(u)\cdot\sigma|^k |f_l(v^*)|^{p/q'}\right)\left(\frac{|u|^\gamma}{\la v^*\ra^l}|g_k(v_1^*)|^{q/p'}\mathds{1}_{|v|>3}\mathds{1}_{|v_1^*|^2>E/2}\mathds{1}_{|\hat{v_1^*}\cdot\sigma|>\sqrt{3}/2}\right)\notag\\
&\hspace{7cm}\times\left(|f_l(v^*)|^{p/r}|g_k(v_1^*)|^{q/r}\right)\,d\sigma\,dv_1\notag\\
&=\int_{\R^3\times\S^2}\left(\frac{|\hat{R}_\sigma^+(u)\cdot\sigma|^k}{|\hat{R}_\sigma^-(u)\cdot\sigma|^{(1-s)\alpha/p'}}|f_l(v^*)|^{p/q'}\right)\left(\frac{|u|^\gamma |\hat{R}_\sigma^-(u)\cdot\sigma|^{\alpha/p'}}{\la v^*\ra^l}|g_k(v_1^*)|^{q/p'}\mathds{1}_{|v|>3}\mathds{1}_{|v_1^*|^2>E/2}\mathds{1}_{|\hat{v_1^*}\cdot\sigma|>\sqrt{3}/2}\right)\notag\\
&\hspace{7cm}\times\left(\frac{|f_l(v^*)|^{p/r}|g_k(v_1^*)|^{q/r}}{|\hat{R}_\sigma^-(u)\cdot\sigma|^{s\alpha/p'}}\right)\,d\sigma\,dv_1\notag\\
&\leq \mathcal{J}_1^{1/q'}\mathcal{J}_2^{1/p'}\left(\int_{\R^3\times\S^2}\frac{|f_l(v^*)|^{p}|g_k(v_1^*)|^{q}}{|\hat{R}_\sigma^-(u)\cdot\sigma|^{s\alpha r/p'}}\,d\sigma\,dv_1\right)^{1/r},\label{Holder nonperp Lr}
\end{align}
where 
\begin{align}
\mathcal{J}_1&=\int_{\R^3\times\S^2}\frac{|\hat{R}_\sigma^+(u)\cdot\sigma|^{kq'}}{|\hat{R}_\sigma^-(u)\cdot\sigma|^{(1-s)\alpha q'/p'}}|f_l(v^*)|^p\,d\sigma\,dv_1\label{I1 nonperp Linfty}\\
\mathcal{J}_2&=\int_{\R^3\times\S^2}\frac{|u|^{\gamma p'}|\hat{R}_\sigma^-(u)\cdot\sigma|^\alpha}{\la v^*\ra^{l p'}}|g_k(v_1^*)|^q\mathds{1}_{|v|>3}\mathds{1}_{|v_1^*|^2>E/2}\mathds{1}_{|v^*|<|v|/6}\,d\sigma\,dv_1.\label{I2 nonperp Linfty}
\end{align}
Since $kq'>1$ and $(1-s)\alpha q'/p'<2$, estimate \eqref{averaging estimate on I_2 Maxwell} implies
\begin{equation}\label{J1 bound Lr}
\mathcal{J}_1\lesssim \|f_l\|_{L^p}^p.
\end{equation}
Since $\alpha>1$ and $l>2/q'$, the same argument as in the estimate for $\mathcal{J}_2$ in the proof of Proposition \ref{Q1 proposition Linf} yields
\begin{equation}\label{J2 bound Lr}
\mathcal{J}_2\lesssim \|g_k\|_{L^q}^q.    
\end{equation}

Then, \eqref{Holder nonperp Lr}, \eqref{J1 bound Lr}, \eqref{J2 bound Lr} imply
$$\la v\ra^k\mathcal{K}_{sin}(f,g)(v)\lesssim \|f\|_{L^p}^{p/q'}\|g\|_{L^q}^{q/p'}\left(\int_{\R^3\times\S^2}\frac{|f(v^*)|^p|g(v_1^*)|^q}{|\hat{R}_\sigma^-(u)\cdot\sigma|^{s\alpha r/p'}}\,d\sigma\,dv_1\right)^{1/r}.$$
Raising the above bound to the $r-$th power, integrating in $v$, and using the exponents' relation, we obtain
\begin{align*}
\|\la v\ra^k\mathcal{K}_{sin}(f,g)\|_{L^r_v}^r&\lesssim \|f\|_{L^p}^{rp/q'}\|g\|_{L^q}^{rq/p'}\int_{\R^3\times\R^3\times\S^2}\frac{|f(v^*)|^p\|g(v_1^*)|^q}{|\hat{R}_\sigma^-(u)\cdot\sigma|^{s\alpha r/p'}}\,d\sigma\,dv_1\,dv\\
&= \|f\|_{L^p}^{r-p}\|g\|_{L^q}^{r-q}\int_{\R^3\times\R^3\times\S^2}\frac{|f(v^*)|^p\|g(v_1^*)|^q}{|\hat{R}_\sigma^-(u)\cdot\sigma|^{s\alpha r/p'}}\,d\sigma\,dv_1\,dv.
\end{align*}
Since $\frac{s\alpha r}{p'}<2$, we use Lemma \ref{lem:singular involutionary estimate} to estimate the remainder term by $\|f\|_{L^p}^p\|g\|_{L^q}^q$, therefore 
\begin{equation}\label{estimate Ksin Lr}
\|\la v\ra^k\mathcal{K}_{sin}(f,g)\|_{L^r_v}\lesssim
\|f\|_{L^p}\|g\|_{L^q},
\end{equation}

 Combining \eqref{estimate Ksm Lr}, \eqref{estimate Kreg Lr}, \eqref{estimate Ksin Lr} \eqref{decomposition of Q}, estimate \eqref{L inf bound Q_1} follows.
The proof is complete.
\end{proof}

\subsection{Proof of Theorem \ref{main theorem}} Let $1\leq p,q<3$ with $p\geq \frac{2}{2-\gamma}$ and $\frac{2}{2-\gamma}\leq r\leq\infty$, such that $\frac{1}{p}+\frac{1}{q}=1+\frac{1}{r}$. Fix $k\geq l>2/p'$ with $k>1/q'$. Then by Propositions \ref{Q1 proposition Linf}, \ref{Q_1 proposition Lr}, we have
\begin{equation}\label{cc bound hard}
\|\la v\ra^k Q^+(f,g)\|_{L^r}\lesssim \|\la v\ra^l f\|_{L^p}\|\la v\ra^kg\|_{L^q},\qquad\forall\, f,g\in C_c(\R^3).    
\end{equation}
The claim then follows by a standard density argument.

\end{document}